\documentclass[12pt,twoside]{amsart}
\usepackage{amssymb, amsmath, mathrsfs, mathtools, amsfonts, amsthm, amscd, yfonts}
\usepackage[a4paper,top=2 cm,bottom=2 cm,left=2cm,right=2cm]{geometry}
\usepackage{xcolor}
\usepackage{cite}
\usepackage{hyperref} 
\usepackage{mathrsfs}
\usepackage{multicol}
\usepackage{mathtools,tikz-cd}
\usepackage[hyperpageref]{backref} 

\usepackage{graphicx,color}
\usepackage{mathrsfs}
\usepackage{pdfpages}
\usepackage{tabularx}

\usepackage{geometry}
\usepackage{color}
\usepackage[latin1]{inputenc}
\usepackage[all]{xy}

\newcommand{\C}{{\mathbb{C}}}
\newcommand{\Z}{{\mathbb{Z}}}
\newcommand{\wP}{{\mathbb{P}}}

\def\C{\mathbb C}
\newcommand{\CC}{{\mathbb C}}

\newcommand{\ddd}{ D}

\renewcommand{\dim}{\mathrm{dim}}

\newcommand{\OO}{\mathcal O}
\newcommand{\F}{\mathscr F}
\newcommand{\Ls}{\mathscr L}
\newcommand{\J}{\mathscr J}
\newcommand{\fol}{\mathscr F}
\newcommand{\G}{\mathscr G}

\renewcommand{\O}{{\mathcal O}}

\newcommand{\PP}{\mathbb{P}}

\newcommand{\reg}{\mathrm{reg}\,}

\newcommand{\singf}{{\rm Sing}({\mathscr F})}

\def\Z{\mathbb Z}

\newcommand{\Sing}{{\rm Sing}}
\newcommand{\Log}{{\rm Log}}
\newcommand{\Ind}{{\rm Ind}}
\newcommand{\Res}{{\rm Res}}
\newcommand{\BB}{{\rm BB}}
\newcommand{\GSV}{{\rm GSV}}
\newcommand{\CS}{{\rm CS}}
\newcommand{\cod}{{\rm cod}}

\newtheorem{lema}{lema}[section]
\newtheorem{cor}[lema]{Corollary}
\newtheorem{teo}[lema]{Theorem}

\theoremstyle{definition}
\newtheorem{remark}[lema]{Remark}

\newtheorem{defi}[lema]{Definition}
\newtheorem{exe}[lema]{Example}

\usepackage{fancyhdr}

\usepackage{lipsum}

\begin{document}

\pagestyle{fancy}


\title[Log Baum--Bott Residues  for foliations by curves]{Log Baum--Bott Residues  for foliations by curves}
\author{Maur\'icio Corr\^ea}
\address{Maur\'icio Corr\^ea  \\ Department of Mathematics,  
Universit\`a degli Studi di Bari, 
Via E. Orabona 4, I-70125, Bari, Italy
}
\email[M. Corr\^ea]{mauricio.barros@uniba.it,mauriciomatufmg@gmail.com }

\author{Fernando Louren\c co}
\address{ Fernando Louren\c co \\ DMM - UFLA ,  Campus Universit\'ario, Lavras MG, Brazil, CEP 37200-000}
\email{fernando.lourenco@ufla.br}

\author{Diogo Machado}
\address{Diogo Machado \\ DMA - UFV,  Avenida Peter Henry Rolfs, s/n - Campus Universitário, 36570-900 Vi\c cosa- MG,
Brazil}
\email{ç}




\begin{abstract}
We prove a Baum-Bott type    residual formula  for  one-dimensional holomorphic foliations, and  logarithmic along   free    divisors.  More precisely, this provides a Baum--Bott theorem for a foliated triple \( (X, \fol, D) \), where \( \fol \) is a foliation by curves and \( D \) is a free divisor on a complex manifold \( X \). From the local point of view, we show that the log Baum--Bott residues are a generalization of the  Aleksandrov logarithmic index for vector fields with isolated singularities on hypersurfaces. We also show how these new indices are related to Poincar\'e's Problem for foliations by curves. In the case of foliated surfaces, we show that the differences between the logarithmic residues and Baum--Bott indices along   invariant curves can be expressed in terms of the GSV and Camacho--Sad indices. We also  obtain  a Baum--Bott type formula for singular varieties via log resolutions. Finally, we prove a weak global version of the Zariski--Lipman conjecture for compact algebraic surfaces, in the form of a foliated smoothness criterion, suggesting the appearance of saddle-nodes in the singularity reduction on singular surfaces. 
\end{abstract}

\maketitle

\hyphenation{sin-gu-lari-ties}
\hyphenation{re-si-dues}



 \section{Introduction}

Poincar\'e, in his foundational memoir \cite{Poincare1885} on the qualitative theory of ordinary differential equations, introduced the concept of the index (or \textit{indices du cycle} at \textit{equilibrium points}) for singular points of planar vector fields. This concept was later further developed and formalized by Hopf \cite{Hopf1927}, leading to the classical  \textit{Poincar\'e-Hopf} Index  Theorem. This theorem established that the sum of the indices of singularities in a vector field on a closed, orientable manifold corresponds to its Euler characteristic, thereby providing a profound link between topology and dynamics.

The extension of the index theorem to the holomorphic setting began with the work of Bott \cite{Bott67}.  He introduced new indices (residues) and demonstrated that the sum of these indices at singularities is equal to the characteristic numbers of the manifold. Building upon Bott's work, Baum and Bott in  \cite{BB1,BB2} further advanced the theory by investigating the residues of singularities in   holomorphic foliations, particularly in the context of meromorphic vector fields with isolated singularities. More precisely, if  $\fol$ is  holomorphic  foliation of dimension one(\textit{foliation by curves})  with  only  isolated singularities on an $n$-dimensional manifold $X$. The  Baum-Bott residues can be expressed in terms of a Grothendieck residue as follows: if  $p \in \Sing(\fol)$ is a singularity  and $\varphi$ a homogeneous symmetric polynomial of degree $n$,   then 
\begin{eqnarray}\label{0001}
\BB_{\varphi}(\fol,  p ) = \mbox{Res}_p\left[\begin{array}{cccc} \varphi(Jv)\\ v_1,\ldots,v_n\end{array}\right ],
\end{eqnarray}
\noindent where $v=(v_1,\ldots,v_n)$ is a germ of holomorphic vector field at $p$, local representative of $\fol$, and $Jv$ is its Jacobian matrix. Moreover, if $X$ is compact, then
$$ \varphi(\mathcal{N_{\fol}}) \cap [X]= \varphi(T_X - T_{\fol})\cap [X] = \sum_{  p \in \Sing(\fol) } \mbox{Res}_p\left[\begin{array}{cccc} \varphi(Jv)\\ v_1,\ldots,v_n\end{array}\right ],$$ 
where $\mathcal{N_{\fol}}$ is the normal sheaf of $\fol$. In particular, if $\fol$ is induced by a global holomorphic vector field and $ \varphi=\det$, then by Chern--Gauss--Bonnet and Poincar\'e--Hopf Thereoms
$$  \chi(X)= c_n(T_X)\cap [X] = \sum_{  p \in \Sing(\fol) } \mbox{PH}_p(\fol),$$
where $ \mbox{PH}_p(\fol)$ is the Poincar\'e--Hopf index of $\fol$ at $p$, and we know that it is equal to
$$
\mbox{PH}_p(\fol)= \mbox{Res}_p\left[\begin{array}{cccc} \det(Jv)\\ v_1,\ldots,v_n\end{array}\right ]= \dim \frac{\mathcal{O}_p}{<v_1,\ldots,v_n>}= \displaystyle\sum_{i=0}^n(-1)^i \dim\  H_i(\Omega^\bullet_{X,p},i_v),
$$
where $H_i(\Omega^\bullet_{X,p},i_v)$   is the homology of the Koszul complex associated with   $v$, which, in turn, also coincides with $\mu_x(\fol)$ the Milnor number of the vector  field $v$ at \( p \).

  Aleksandrov \cite{Ale_p}, using homological algebra, has defined a certain algebraic index which measures the variation between the Poincar\'e-Hopf and homological indices as follows: Let $v$ be a local vector field inducing $\fol$. The interior multiplication $i_{v}$ induces the complex of logarithmic differential forms
$$
0 \longrightarrow \Omega^n_{X,p}(\log\, \ddd) \stackrel{i_v}{\longrightarrow} \Omega^{n-1}_{X,p}(\log\, \ddd)\stackrel{i_v}{\longrightarrow} \cdots\stackrel{i_v}{\longrightarrow} \Omega^{1}_{X,p}(\log\, \ddd)\stackrel{i_v}{\longrightarrow} \OO_{n,p} \longrightarrow 0,
$$
where $D$ is a divisor on $X$. If $p$ is an isolated singularity of $\fol$, the $i_v$-homology groups of the complex $\Omega^{\bullet}_{X,p}(\log\, \ddd)$ are finite-dimensional vector spaces (see \cite{Ale_p}). Thus, the Euler characteristic
$$
\chi (\Omega^{\bullet}_{X,p}(\log\, \ddd), i_v) = \displaystyle\sum_{i=0}^n(-1)^i \dim H_i(\Omega^\bullet_{X,p} (\log\, \ddd),i_v)
$$
\noindent of the complex of logarithmic differential forms is well defined.  Since this number does not depend on local representative $v$ of the foliation $\fol$ at $p$, the {\it logarithmic index} of $\fol$  at the point $p$ is defined by
$$
\Log(\fol,\ddd,p):= \chi (\Omega^{\bullet}_{X,p}(\log\, \ddd), i_v).
$$
\noindent  Aleksandrov \cite{Ale_p} has proved the following formula that relates the logarithmic index and the residues  
\begin{eqnarray} \label{0003}
\Log(\fol, \ddd, p) = \mbox{Res}_p \left[\begin{array}{cccc} \det(Jv) \\ v_1,\ldots,v_n \end{array}\right] - \mbox{Ind}_{Hom} (v,  D, p),
\end{eqnarray}
where $\mbox{Ind}_{Hom}(v,  D, p)$ is the so-called  homological index of $v$ at $p$.

Let $X$ be a $n$-dimensional complex manifold, $\ddd$ an isolated hypersurface singularity on $X$ and let $\fol$ be a foliation on $X$ of dimension one, with isolated singularities. Suppose  that $\fol$ is  logarithmic along $\ddd$, i.e., the analytic hypersurface $\ddd$ is invariant by each holomorphic vector field that is a local representative of $\fol$. The $\GSV$-index of $\fol$ in $x\in \ddd$ will be denoted by $\GSV(\fol,\ddd,x)$. Recall that $\mbox{Ind}_{Hom}(v,  D, x)=\GSV(v,D,x)$ if $x$ is isolated,  see \cite{GM1,GSV}. 
For definition and details on the $\GSV$-index we refer to \cite{bss} and \cite{Suw2}. See also Subsection \ref{GSV-CS}.

A complex $\partial$-manifold is a complex manifold of the form $\tilde{X} = X - D$, where $X$ is an $n$-dimensional complex compact manifold and $D \subset X$ is a divisor referred to as the boundary divisor. Nowadays, the following version of the   Gauss-Bonnet theorem for $\partial$-manifolds is well known\cite{Aluffi}:
$$
\int_{X} c_n(T_{X}(-\log\, D)) \,\,\, = (-1)^n\,\,\, \chi(\tilde{X})=\displaystyle \sum_{i=1}^{n} \dim H^i_{c}(\tilde{X}, \CC).
$$

In \cite{CM1}, the first and third authors addressed the problem of providing a \textit{log Poincar\'e--Hopf} theorem for meromorphic vector fields on compact complex manifolds and proved the following result under the assumption that invariant divisors are of normal crossing type.

\begin{teo}\label{teo_1.1}
Let $\fol$ be a one-dimensional foliation with isolated singularities and logarithmic along the normal crossing divisor $D$ in a complex compact manifold $X$.  Then
$$
\displaystyle\int_{X}c_{n}(T_{X}(-\log\, D)- T_{\fol}) =    \sum_{x\in  \Sing(\fol)\cap (X\setminus  D)} \mu_x(\fol) + \sum_{x\in \Sing(\fol) \cap  D }   \Log(\F,D,p),
$$
where $T_{\fol}$ denotes  the tangent bundle of $\fol$ and  $ \mu_x(\fol)$ is the Milnor number of  $\fol$ at $x$. 
\end{teo}

The main goal of this article is to establish a Baum-Bott type theorem for foliations by curves on pairs \((X, D)\) and  more general characteristic numbers of foliations, especially when \( D \) is a free divisor. The two primary motivations are the following: to extend the Baum-Bott theorem to non-compact $\partial$-manifolds, and to singular normal projective algebraic varieties via resolutions of singularities. 
  This type of problem was first considered by Seade and Suwa in \cite{SeaSuw1}. They presented residue formulas for foliations in local complete intersection varieties with isolated singularities (ICIS).
   Also, a  Bott-type formula for complex orbifolds is studied in \cite{CorreaRodriguezSoares2016}, and formulas relating residues via good resolutions are obtained.

  We propose a more general approach by using techniques that take into account the pair \( (X, D) \) and the respective logarithmic tangent sheaf, allowing for the inclusion of singularities more general than ICIS.

 Let $\fol$ be a one-dimensional holomorphic foliation on a complex manifold $X$, logarithmic along an analytic free divisor $\ddd \subset X$. Given $p \in \Sing(\fol) \cap \ddd$ an isolated singular point of $\fol$, let $U\subset X$ be a neighborhood of $p$ and $\vartheta \in T_X(-\log \,\ddd)|_{U}$ a local representative of $\fol$. Setting $(z_1,\ldots,z_n)$ a system of complex coordinates on $U$, we can write the vector field $\vartheta$ explicitly as 
\begin{eqnarray}
\vartheta = \displaystyle \sum_{j=1}^n v_j \frac{\partial}{\partial z_j}.
\end{eqnarray}

\noindent  
Let $\delta_1,\dots,\delta_n$ be a local $\mathcal O_X$-frame of $T_X(-\log D)$ on $U$,
and write
\[
[\delta_k,\delta_j]=\sum_{i=1}^n c^i_{kj}\,\delta_i,\qquad
\vartheta=\sum_{k=1}^n \vartheta^k \delta_k.
\]
Define the matrix $M_{\log}(\vartheta)=(m_{ij})_{1\le i,j\le n}$ by
\begin{equation}\label{eq:Mlog-clean}
m_{ij}:=-\delta_j(\vartheta^i)+\sum_{k=1}^n \vartheta^k\,c^i_{kj}.
\end{equation}
Then
\[
[\vartheta,\delta_j]=\sum_{i=1}^n m_{ij}\,\delta_i,
\]
so that $M_{\log}(\vartheta)$ is the matrix of the action $\alpha_\vartheta$ in the frame
$\delta_1,\dots,\delta_n$.
If $D$ has normal crossings on $U$ and the chosen frame is the standard logarithmic
frame, then $[\delta_k,\delta_j]=0$ and \eqref{eq:Mlog-clean} reduces to
$M_{\log}(\vartheta)=J_{\log}\vartheta$.

Let $\varphi$ be a homogeneous symmetric polynomial of degree $n$ and $\fol$ a one-dimensional holomorphic foliation with isolated singularities.   We generalize the logarithmic index by defining, for each $p \in \Sing(\fol) \cap \ddd$, the \textit{Log Baum--Bott index} of $\fol$ at $p$ with respect to $\varphi$ as follows: 
$$\mbox{Res}^{\log }_{\varphi}(\fol,D, p) := \mbox{Res}_p\displaystyle\left[\begin{array}{cccc}\varphi(M_{\log} (\vartheta))dz_{1}\wedge \cdots \wedge dz_{n} \\ v_1,\ldots,v_n \end{array}\right ],$$
where $\mbox{Res}_p$ means the  Grothendieck residue simbol at $p$ with respect $v_1,\ldots,v_n$.

We prove the following Baum--Bott residual type theorem for foliations by curves that are logarithmic along an analytic free divisor.

\begin{teo} \label{ttt1}  Let $\fol$ be a one-dimensional holomorphic foliation on a complex manifold $X$, logarithmic along a  free divisor $\ddd \subset X$ and $\varphi$ a homogeneous symmetric polynomial of degree $n$. Then:
\begin{itemize}
\item [(i)]  for each connected component $S_{\lambda} \subset \Sing(\fol)\cap D$ of the singular set of $\fol$, there exists the residue $\Res^{\log}_{\varphi}\Big(\fol,D,  S_{\lambda}\Big)$,  which is a complex number that depends only on $\varphi$ and the local behavior of the leaves of $\fol$ near $S_{\lambda}$.
\item [(ii)] If $X$ is compact, then
$$\int_{X} \varphi\Big(T_X(-\log\, \ddd) - T_{\fol}\Big) = \sum_{ S_{\lambda} \not\subset D } \BB_{\varphi}\Big(\fol, S_{\lambda}\Big)+ \sum_{S_{\lambda}\subset  D   } \Res^{\log}_{\varphi}\Big(\fol,   \ddd  , S_{\lambda}\Big),$$
\noindent where the sum is taken over all the connected components $S_{\lambda}$ of the singular set of $\fol$ and $\BB_{\varphi}\Big(\fol, S_{\lambda}\Big)$ denotes the Baum-Bott residue of $\fol$ at $S_{\lambda}$ with respect to $\varphi$.
\item [(iii)] Given $p \in \Sing(\fol) \cap \ddd$ an isolated singular point of $\fol$,  let $\vartheta$ be a logarithmic vector field that is a local representative of $\fol$ (in a neighborhood of $p$). Then
$$\Res^{\log }_{\varphi}( \fol,D, p) = \Res_p\displaystyle\left[\begin{array}{cccc}\varphi\big(M_{\log} (\vartheta)\big)  \\ v_1,\ldots,v_n \end{array}\right ],$$
\noindent where $M_{\log}(\vartheta)$ denotes the $n\times n$ matrix associated with $\vartheta$ defined by $(\ref{eq:Mlog-clean})$  and $v_1,\ldots,v_n$ denote the local coordinates of
$\vartheta$ as a section of $T_X$. 
\end{itemize}
\end{teo}
 
\begin{remark} For the index, residues, and other related topic, we use the notations appearing in Suwa's book in \cite{Suw2}.
    
\end{remark}

From Theorem \ref{ttt1} we conclude the following.

\begin{cor} \label{cor-princ-iso}  If $\fol$ is a foliation by curves  with only  isolated singularities, then  
$$
\displaystyle\int_{X}\varphi(T_{X}(-\log\, D)- T_{\fol}) =    \sum_{x\in  \Sing(\fol)\cap (X\setminus  D)}  \BB_{\varphi}(\fol, x)+ \sum_{x\in \Sing(\fol) \cap  D }   \Res^{\log }_{\varphi}(\fol, p),
$$
where   $ \BB_{\varphi}(\fol, x)$ is the Baum-Bott residue   of  $\fol$ at $x$ with respect to $\varphi$. 
If in addition, $\ddd$ is  normal crossing at $p$,  
then
$$\Res^{\log }_{\varphi}( \fol,D, p) = \Res_p\displaystyle\left[\begin{array}{cccc}\varphi(J_{\log} \vartheta)  \\ v_1,\ldots,v_n \end{array}\right ].$$
 Also, when $\varphi = \det$ we have 
$$\mbox{\Res}^{\log }_{\det}(\fol,D, p) =  \Log(\F,D,p).$$
\end{cor}

Let $W\subset X$ be an analytic
subspace and $\fol$ a holomorphic foliation with isolated singularities, we denote by
$$
\mbox{R}(\fol,W)=\sum_{x\in \Sing(\fol)\cap W} \mbox{R}(\fol,x),
$$
where $\mbox{R}$ denotes some residue associated with $\fol$ along $W$. If \( \mathrm{R} = \Res^{\log}_{\phi} \), for some homogeneous symmetric polynomial \( \phi \), and \( W \) is not contained in \( D \), we still use \( \Res_{\phi}(\fol , W) = \Res^{\log}_{\phi}(\fol, W) \), where \( \Res_{\phi}(\fol, W) \) is the total sum of the Baum--Bott residues of \( \fol \) along \( W \).

In \cite{Poincare}, H. Poincar\'e observed that establishing an upper bound on the degree of algebraic solutions would be sufficient to ensure the existence of local first integrals for polynomial ordinary differential equations.
 This question, now referred to as Poincar\'e's Problem, plays a fundamental role in the study of  holomorphic foliations.  A positive answer to Poincar\'e's Problem can be obtained by imposing some  conditions on either the foliation or the invariant variety. We refer the reader to \cite{Correa2024} and the references therein.

We found a new obstruction to giving an affirmative answer to the Poincar\'e problem for foliations by curves on $\mathbb{P}^n$, with $n$ odd. 
\begin{cor}\label{Poincare-oddn}
    Let $\fol$ be a one-dimensional foliation, with isolated singularities, on $\mathbb{P}^n$ with an invariant free divisor $D$. If  $n$ is odd   and 
    $$
 \Res^{\log }_{c_{1}^n}(\fol, \Sing(\fol))\geq 0
$$
then 
\begin{eqnarray}\label{poincnodd}
\deg(D) \leq  \deg(\fol) +n.
\end{eqnarray}
\end{cor}
The bound (\ref{poincnodd}) was obtained, under mild conditions, by Soares \cite{Soares}, Esteves \cite{Esteves}, and Brunella and Mendes \cite{Brunella-Mendes}.

Now,   let $\fol$ be a one-dimensional holomorphic foliation on a compact complex surface $X$, logarithmic along a reduced curve $\ddd \subset X$. Then, by corollary \ref{cor-princ-iso} we have 
   $$
 c_2(T_{X}(-\log\, D)- T_{\fol})\cap [X]=  \mu(\fol, X\setminus D)+ \Res^{\log }_{\det}( \fol,D)
$$
and 
  $$
 (c_1(N_\fol)-D)^2= \BB_{c_1^2}(\fol,   X\setminus D)+  \BB_{c_1^2}(\fol,   D ).$$
 
In the following result, we will compare the differences between the logarithmic indices and the Baum--Bott indices along \( D \) in terms of the GSV
and Camacho--Sad indices. We  recall that  a foliation \( \fol\) on a complex surface $X$ is a \textit{generalized curve} along an invariant curve \( D \) if, for all \( p \in \Sing(\fol) \), \( p \) is non-dicritical and there are no saddle-nodes in its resolution.

\begin{cor}\label{dim2-poinc}
 Let $\fol$ be a one-dimensional holomorphic foliation on a compact complex surface $X$, logarithmic along a reduced curve $\ddd \subset X$. Then
$$
\BB_{c_1^2}(\fol, D)-\Res^{\log }_{c_1^2}(\fol, D) =2\GSV(\fol,D)+\CS(\fol,D).
$$
In particular,  the  non-negativity  of 
$$
\BB_{c_1^2}(\fol, D)-\Res^{\log }_{c_1^2}(\fol, D) -\CS(\fol,D)
$$
is an  obstruction to an affirmative answer to the Poincar\'e problem.  Moreover, if $\fol$  is a generalized curve along $D$, then
$$
\Res^{\log }_{c_1^2}(\fol, D)=0.$$
\end{cor}

\noindent This formula can be illustrated locally as follows: Consider \( v = \lambda z\frac{\partial}{\partial z} + \mu w\frac{\partial}{\partial w} \), with \( D = \{ z = 0 \} \) such that \( \lambda, \mu \neq 0 \). Then
\[
\Res^{\log}_{c_1^2}(v, D, 0) = \frac{\operatorname{tr}^2(J_{\log} v)}{\det(Jv)} = \frac{\mu}{\lambda} = \CS(D, 0)^{-1}
\]
and
\[
\BB_{c_1^2}(v, 0) - \Res^{\log}_{c_1^2}(v, D, 0) = \frac{\lambda}{\mu} + 2 = \CS(v, D, 0) + 2.
\]
So,
\[
\BB_{c_1^2}(v, 0) - \Res^{\log}_{c_1^2}(v, D, 0) - \CS(v, D, x) = 2 = 2 \, \GSV(v, D, x).
\]

In \cite{SeaSuw1}, Seade and Suwa derive Baum--Bott type residue formulas for foliations by curves induced by global holomorphic vector fields on complex manifolds with nonempty boundary. In the orbifold setting, the  first author, Rodr\'iguez and Soares obtain an analogous Bott-type residue formula and, moreover, compare local residue contributions through suitable (good) resolutions \cite{CorreaRodriguezSoares2016}. In the sequel we prove a more general statement in our logarithmic framework, using log resolutions to relate global characteristic numbers to sums of logarithmic Grothendieck--Baum--Bott residues, as follows.

  As usual, we denote by $T_Y:= \mbox{Hom}( \Omega_Y^1, \O_Y)$ the tangent sheaf of $Y$, where $\Omega_Y^1$ is the sheaf of   K\"ahler differentials of $Y$. 
  
  We obtain the following result.
\begin{cor} \label{Baum-Bott-sing}
Let $Y$ be a compact normal variety,   $ T\fol \subset T_Y$  a one-dimensional holomorphic foliation,  with isolated singularities, and $\varphi$ a homogeneous symmetric polynomial of degree $n=\dim(Y)$.
If $\pi:(X,D)\to (Y, \emptyset)$  is  a functorial   log  resolution of $Y$, then 
 $$
 \int_{X}\varphi(T_{X}(-\log\, D)-(T\pi^*\fol)^*) =     \sum_{x\in  \Sing(v)\cap Y_{\reg} }    \BB_{\varphi}(v, x) +  \sum_{ S_{\lambda} \subset D\cap \Sing(\pi^*\fol)} \Res^{\log}_{\varphi} (\pi^*\fol, S_{\lambda} ),
  $$
  where $ \BB_{\varphi}(\pi^*\fol, x)$ is the Baum--Bott index  of  $\pi^*\fol$ at $x$ with respect to $\varphi$. In particular,   if   $ D\cap \Sing(\pi^{-1}v)$ is isolated, then

  $$
 \int_{X}\varphi(T_{X}(-\log\, D)- (T\pi^*\fol)^*) =     \sum_{x\in  \Sing(\fol)\cap Y_{reg} }    \BB_{\varphi}(\fol, x)+  \sum_{ x\in  D\cap \Sing(\pi^{-1}v)}  \Res_p\displaystyle\left[\begin{array}{cccc}\varphi(J_{\log}(\pi^{-1}v))  \\ v_1,\ldots,v_n \end{array}\right ].
  $$

  If, in addition, \( \pi^*\fol \) is induced by a global vector field $v$ and \( \Sing(\fol) \subset \Sing(Y) \), we obtain a \textit{Poincar\'e--Hopf} type Theorem
 $$
  \chi (Y)=  \sum_{ x \in D\cap \Sing(\pi^{-1}v)} \Res_p\displaystyle\left[\begin{array}{cccc}\varphi(J_{\log}(\pi^{-1}v))  \\ v_1,\ldots,v_n \end{array}\right ]+  \chi (\Sing(Y)).
 $$
 In particular, for each $\varphi$,  the number  
\[
\sum_{x \in D \cap \Sing(\pi^{-1}v)} \Res_p\displaystyle\left[\begin{array}{cccc}\varphi(J_{\log}(\pi^{-1}v))  \\ v_1,\ldots,v_n \end{array}\right ]
\]
is a topological invariant; it depends only on \(v\) and does not depend on the resolution.
\end{cor}
Observe also  that 
$$ 
\BB_{\varphi}(\fol, Y_{reg}) := \sum_{x \in \Sing(\fol) \cap Y_{\reg}} \BB_{\varphi}(\fol, x) 
$$ 
does not depend on the log resolution \( \pi : (X, D) \to Y \).

Now, let $\pi:(\widehat X,\widehat D)\to (X,D)$ be a  log resolution of the log pair $(X,D)$ with $X$ smooth,  
$\widehat D=\pi^{-1}(D)_{\mathrm{red}}$ is simple normal crossings. Let $\widehat Y:=\widehat D_{\mathrm{sing,sch}}$ be the Jacobian (singular) scheme of $\widehat D$, and set
$L:=\mathcal O_{\widehat X}(\widehat D)$.
Let $v\in H^0\bigl(X,T_X(-\log D)\bigr)$ be a global logarithmic vector field with isolated zeros, and suppose that 
$\widehat v\in H^0\bigl(\widehat X,T_{\widehat X}(-\log \widehat D)\bigr)$ has also isolated zeros.
Consider the $\mu$-class  $\mu_L(\widehat Y)$ with respect to line bundle $L=\mathcal O_{\widehat X}(\widehat D)$ (in the sense of Aluffi, in \cite{Aluffi1}), namely
$$\mu_L(\widehat Y):=c(T^*_{\widehat X}\otimes L)\cap s(\widehat Y,\widehat X).$$

Under the assumptions above, we thus obtain a new theorem of Hopf--Gauss--Bonnet type.

\begin{cor}\label{mu-PHGB}
Then
\begin{equation}\label{eq:log-index-mu-class-short}
\sum_{p\in \Sing(\widehat\fol)\cap \widehat D} \Log(\widehat\fol,\widehat D,p)
=(-1)^{n}
 \int_{\widehat X}\mu_L(\widehat Y).
\end{equation}
\end{cor}

The Zariski--Lipman conjecture   \cite{Lip}  predicts  that an analytic variety with a locally free tangent sheaf is smooth. This conjecture holds for varieties with singularities of codimension $\geq 3$ \cite{Fle88}; local complete intersection singularities \cite{Kal11},  klt   singularities \cite{Greb-Keb-Kov-Pet-11,AD}; log canonical singularities \cite{Greb-Keb-Kov-Pet-11,Dru14, GraKov14}.
 Furthermore, conjecture is established for certain types of surfaces \cite{Bal06, RicOel88, Graf1,GraBer20, BisGurKol14}.  In \cite{Graf1}, Graf proposed the study of a global version of the Zariski--Lipman conjecture and proved that the conjecture holds under certain geometric conditions on the surface. In \cite{GraBer20}, Bergner and Graf prove that a global version of the conjecture is also true if the surface is compact and has two twisted vector fields that are linearly independent at some point. It is natural to suppose that a compact surface \( X \), with a locally free tangent sheaf, has a global twisted vector field that is non-singular on its regular part; that is, it is at most singular on the singular set of \( X \), see \cite[page 2]{Graf1}. We prove the following weak global version.

Consider    $\pi:(X,\tilde{D})\to Y$  a    log  resolution of $Y$
 and a twisted vector
field $v\in H^0(Y, T_Y\otimes \Ls)$ such that $\Sing(v)\subset \Sing(Y)$. If  $\G=\pi^*\fol$  is the induced foliation on $X$, then $\Sing(\pi^*\fol)\subset D$.  Now, consider $\rho:(Z,D)\to X$ a resolution of singularities of $\G$ . So, up to take a singularity resolution of the foliation, we can suppose that all singularities of $\G$ are reduced, i.e., if $p\in \Sing(  \G)$, then $p$ is non-degenerate or is a saddle-node. 
Recall that  $\fol$  is a generalized curve if the resolution $\G$ does not contain saddle-nodes \cite{Bru97}. 
We obtain the following result.
\begin{cor}\label{GZL}
   Let $Y$ be a compact analytic surface such that and $\fol$ a foliation by curves on $Y$ induced by a twisted vector field $v \in H^0(Y, T_Y \otimes \mathscr{L})$, where $\mathscr{L}$ is an invertible sheaf and $\Sing(\fol) \subset \Sing(Y)$. If $TY$ is locally free and  \(\fol \) is a generalized curve,   then $Y$ is smooth. 
\end{cor}

This yields a foliated smoothness criterion: if $v\in H^0(Y,\,T_Y\otimes\mathscr L)$ is a twisted vector field and the induced foliation $\fol$ satisfies $\Sing(\fol)\subset \Sing(Y)$, then $\fol$ cannot be a generalized curve unless $Y$
is smooth. Equivalently, if $Y$ is singular while $T_Y$ is locally free, then for every such $v$ the reduction of singularities of
$\fol$ necessarily exhibits at least one saddle-node.
 In this sense, the statement may be viewed as a foliated Zariski--Lipman phenomenon: the algebro-geometric condition that
$T_Y$ be locally free manifests itself as a rigid dynamical constraint on foliations arising from twisted vector fields
whose singular locus is confined to $\Sing(Y)$.
 
In a recent work, Liao and Zhang \cite{LiaoZhang} give a microlocal interpretation of several local indices for logarithmic
vector fields, defining them as local intersection numbers between characteristic cycles of constructible functions and the
graph of the covector field associated  to a logarithmic vector field. This framework is powerful
enough to recover, and relate,   GSV, and logarithmic indices through relations among standard constructible
functions. They also stress an intrinsic
limitation of this purely intersection-theoretic mechanism: because characteristic cycles carry integral multiplicities,
the resulting indices are necessarily integer-valued, so complex-valued invariants of Camacho--Sad type cannot be captured
in this way. Motivated by this observation, in the final section \ref{sec:residue-constructible} we explain how our logarithmic
Grothendieck--Baum--Bott residues   can still be organized by a natural residue-valued constructible
function.

\subsection*{Acknowledgments} We would   like to thank Alan Muniz and Jos\'e Seade for useful discussions. 
MC is partially supported by the Universit\`a degli Studi di Bari and by the
 PRIN 2022MWPMAB- "Interactions between Geometric Structures and Function Theories" and he is a member of INdAM-GNSAGA;
he was  partially supported by CNPq grant numbers 202374/2018-1, 400821/2016-8 and  Fapemig grant numbers APQ-02674-21,  APQ-00798-18,  APQ-00056-20. FL  is  partially supported by Fapemig grant number APQ-02674-21.

\section{Preliminaries}
\subsection{Logarithmic forms and logarithmic vector fields}

Let $X$ an $n$-dimensional complex manifold and $\ddd$ a reduced hypersurface on $X$. Given a  meromorphic $q$-form $\omega$ on $X$, we say that $\omega$ is a {\it logarithmic} $q$-form along $\ddd$ at $p\in X$ if the following conditions occurs:

\begin{enumerate}
\item[(i)] $\omega$ is holomorphic on $X - \ddd$;
\medskip
\item[(ii)] If $h=0$ is a reduced equation of $\ddd$, locally at $p$, then $h \,\omega$ and $h\, d \omega$ are holomorphic.
\end{enumerate}
\noindent Denoting by $\Omega_{X,p}^q(\log\, \ddd)$ the set of germs of logarithmic $q$-forms  along $\ddd$ at $x$, we define the following coherent sheaf of $\OO_X$-modules
$$
\Omega_{X}^q(\log\, \ddd):= \bigcup_{p\in X} \Omega_{X,p}^q(\log\, \ddd),
$$
which is called  by {\it sheaf of logarithmic $q$-forms along $\ddd$}.  See \cite{Deligne}, \cite{Katz} and  \cite{Sai}  for details. 

Now, given $p\in X$, let $v\in T_{X,p}$ be germ at $p$ of a holomorphic vector field on $X$. We say that $v$ is a {\it logarithmic vector field along of $\ddd$ at $p$}, if $v$ satisfies the following condition: if $h=0$ is an equation of $\ddd$, locally at $p$, then the derivation $v(h)$ belongs to the ideal $\langle h_p\rangle\OO_{X,p}$. Denoting by $T_{X,p}(-\log\, \ddd)$ the set of germs of logarithmic vector fields along of $\ddd$ at $p$, we define the following coherent sheaf of $\OO_X$-modules
$$
T_{X}(-\log\, \ddd):= \bigcup_{p\in X}T_{X,p}(-\log\, \ddd),
$$
which is called by  {\it sheaf of logarithmic vector fields along $\ddd$}. 

It is known that $\Omega_{X}^1(\log\, \ddd)$ and $T_{X}(-\log \,\ddd)$ are always reflexive sheaves, see  \cite{Sai} for more details. If  $\ddd$ is an analytic hypersurface with normal crossing singularities, the sheaves $\Omega_{X}^1(\log\, \ddd)$ and $T_{X}(-\log \,\ddd)$ are  locally free, furthermore, the Poincar\'e residue map
$$
\Res: \Omega_X^1(\log\, \ddd) \longrightarrow \OO_{\ddd} \cong \bigoplus_{i=1}^N\OO_{\ddd_i}
$$
give the following exact sequence of sheaves on X
\begin{eqnarray}\label{seq1}
0 \longrightarrow \Omega_X^1 \longrightarrow \Omega_X^1(\log\, \ddd) \stackrel{\mbox{Res}}{\longrightarrow} \bigoplus_{i=1}^N\OO_{\ddd_i} \longrightarrow 0,
\end{eqnarray}
\noindent where $\Omega_X^1$ is the sheaf of holomorphic $1$-forms on $X$ and $\ddd_1,\ldots,\ddd_N$ are the irreducible components of $\ddd$.

Let $i:D\hookrightarrow X$ be the inclusion. Define the Jacobian ideal sheaf
$
J_D := \mathrm{Fitt}_1(\Omega_D^1)\subset \mathcal O_D,
$
and view the twist
$
J_D(D):= i_*\bigl(J_D\otimes_{\mathcal O_D}\mathcal O_D(D)\bigr)$
as a coherent sheaf on $X$ supported on $D$.
Locally, if $D=\{f=0\}$ with $f$ reduced, there is a natural morphism
\[
T_X \longrightarrow i_*\mathcal O_D(D),\qquad v \longmapsto \frac{v(f)}{f}\Big|_D,
\]
whose kernel is $T_X(-\log D)$. Its image is precisely $J_D(D)$.
Therefore there is a short exact sequence of sheaves on $X$:
\begin{equation}\label{eq:TX-logD-JD}
0 \longrightarrow T_X(-\log D)\longrightarrow T_X\longrightarrow J_D(D)\longrightarrow 0.
\end{equation}

On the projective space $\PP^n$, if $\ddd$ is a smooth   hypersurface,   then there exists the following exact sequence \cite{EA}:
\begin{eqnarray}\label{p2.07}
0 \longrightarrow T_{\PP^n}(-\log\, \ddd) \longrightarrow \OO_{\PP^n}(1)^{n+1} \longrightarrow \OO_{\PP^n}(k)\longrightarrow 0,
\end{eqnarray}
\noindent where $ k$ is the degree of $\ddd$.

\subsection{Singular one-dimensional holomorphic foliations } \label{s0003}

A \emph{singular holomorphic  foliation} $\fol$ in $X$, of dimension $k$, is a reflexive subsheaf $T_\fol$, of rank $k$, of the tangent sheaf $TX$ of $X$ which is involutive, that is, $[T_\fol, T_\fol] \subset T_\fol$. We have   a short exact sequence
\begin{equation}\label{sequ-def-fol} 
    0 \longrightarrow T_\fol \longrightarrow TX\longrightarrow \mathcal{N}_\fol  \longrightarrow 0,
\end{equation}
where $\mathcal{N}_\fol $ is   called the \emph{normal sheaf} of $\fol$, and $ T_\fol$  is  the \emph{tangent sheaf} of $\fol$. Since $T_\fol$ is a reflexive subsheaf of $TX$, then $\mathcal{N}_\fol$ is torsion-free. 

A one-dimensional foliation $\fol$ is also called \textit{foliation by curves} and $K_\fol=T_{\fol}^*$ is its canonical bundle of $\fol$. Twisting the inclusion $T_\fol \longrightarrow TX$ by $K_\fol$, we obtain a global section $v_{\F}\in \mathrm{H}^0(X,T_X\otimes K_{\F})$, which is equivalent to that 
given by the following data:
\begin{itemize}
  \item[$i)$] an open covering $\mathcal{U}=\{U_{\alpha}\}$ of $X$;
  \item [$ii)$] for each $U_{\alpha}$ a holomorphic vector field $v_\alpha$;
  \item [$iii)$]for every non-empty intersection, $U_{\alpha}\cap U_{\beta} \neq \emptyset $, a
        holomorphic function $$f_{\alpha\beta} \in \mathcal{O}_X^*(U_\alpha\cap U_\beta)$$
\end{itemize}
such that $v_\alpha = f_{\alpha\beta}v_\beta$ in $U_\alpha\cap U_\beta$ and $f_{\alpha\beta}f_{\beta\gamma} = f_{\alpha\gamma}$ in $U_\alpha\cap U_\beta\cap U_\gamma$. 
So $\{f_{\alpha\beta}\}\in \mathrm{H}^1(X, \mathcal{O}^*)$ is a cocycle for $K_{\F}$.  The singular set of $\F$ 
is $\singf=\{v_{\F}=0\}$. The normal sheaf $N_\fol$ is torsion-free, then $\cod(\singf)\geq 2$, and the dual of the exact sequence (\ref{sequ-def-fol}) yields 
$$
 0 \longrightarrow \mathcal{N}_\fol^* \longrightarrow \Omega_X^1\longrightarrow  \mathscr{I}_{\singf}\otimes  K_{\F}  \longrightarrow 0,
$$
where $ \mathscr{I}_{\singf}$ is the sheaf of ideals of $\singf$. If $X$ is a complex surface, we consider the line bundle \( N_{\fol} = (\mathcal{N}_{\fol}^{*})^{*} \), which we will refer to as the normal bundle associated with $\fol$, as usual.

\begin{defi} 
Let $V $ be an analytic subspace of a complex manifold $X$.  We say that $V$ is invariant by a foliation $\F$ if $T_{\F}|_{V}\subset (\Omega^1_V)^*$.
If $V$ is a hypersurface we say that $\F$ is \textit{logarithmic along $V$}.
\end{defi}

Let $\fol$ be a logarithmic one-dimensional holomorphic foliation along a hypersurface $\ddd \subset X$, then
$ T_{\fol} $ is a subsheaf of $ T_X(-\log \,\ddd)$ and $v_{\F}$ is a section of $ \mathrm{H}^0(X, T_X(-\log \,\ddd)\otimes K_{\F})$.

\begin{defi} 
A one-dimensional  foliation on a complex projective space $\mathbb{P}^n$ is called a \textit{ projective foliation}.
Let $\F $ be a projective foliation with tangent bundle $T_{\F}=\mathcal{O}_{\mathbb{P}^n}(r)$. The integer number $d:=r+1$ is called the \textit{ degree} of $\F$.
\end{defi}

Now we can see that Corollary \ref{Poincare-oddn} follows straightforwardly from Theorem \ref{ttt1}. 
\subsection{Proof of Corollary \ref{Poincare-oddn}}
In fact, if $n$ is odd and 
$
 \Res^{\log }_{c_{1}^n}(\fol, \Sing(\fol))\geq 0
$
 then
$$
\deg(D) \leq  \deg(\fol) +n
$$
since  
$$
 (\deg(D) - \deg(\fol) -n)^n =\displaystyle\int_{\mathbb{P}^n}c_1^n(T_{X}(-\log\, D)- T_{\fol}) =    \Res^{\log }_{c_{1}^n}(\fol, \Sing(\fol)).
$$
See the Example \ref{example-arrang}.

\subsection{ Logarithmic,  homological and  $\GSV$ indices  }
Let $D\subset X$  be a reduced hypersurface with  
a local equation $f\in \mathcal{O}_{X,x}$
in a neighborhood of a point $x\in D$.
Consider the $\mathcal{O}_{D,x}$-module of
germs of regular differentials of order $i$ on $D$ :
$$
\Omega_{D,x}^i=\frac{\Omega_{X,x}^i}{f\Omega_{X,x}^i+
df\wedge\Omega_{X,x}^{i-1}}.
$$
Let  $\fol $ be a one-dimensional holomorphic foliation on $X$, with isolated
singularities, logarithmic along $D$. Let $x\in\singf$ be and consider a germ of vector field $v\in T_{X}(-\log \,D)|_{U}$ on $(U,x)$ tangent to $\fol $, where $U$ is a neighborhood of $x$.   Since $v$ is also  tangent to $(D,x)$ the interior multiplication $i_{v}$ induces the complex
$$
0\longrightarrow \Omega_{D,x}^{n-1} \stackrel{i_v}{\longrightarrow} \Omega_{D,x}^{n-2} \stackrel{i_v}{\longrightarrow}\cdots \stackrel{i_v}{\longrightarrow}\Omega_{D,x}^1 \stackrel{i_v}{\longrightarrow}\OO_{D,x}\longrightarrow0.
$$ 
The \emph{homological index} is defined as the Euler characteristic of the complex $(\Omega_{D,x}^{\bullet}, i_v)$
$$
\mbox{Ind}_{Hom}(\fol,D,x)=\sum_{i=0}^{n-1}(-1)^i\dim H_i(\Omega_{D,x}^{\bullet}, i_v). 
$$
Since  the vector field $v$ has an isolated singularity at  $x$, then the $i_v$-homology groups
of the complex $\Omega_{D,x}^{\bullet}$ are finite-dimensional vector spaces and   the Euler characteristic is well defined. Similarly, the Euler characteristic of the Koszul complex 
  $(\Omega_{X,x}^{\bullet}, i_v)$ associated to $v$  is well defined, and the \emph{ Milnor number}  of $\fol$ at $x$ is defined  by 
$$
\mu(\fol,x)=\sum_{i=0}^{n}(-1)^i\dim H_i(\Omega_{X,x}^{\bullet}, i_v). 
$$
 which  coincides with the  Poincar\'e Hopf   index of $v$ at $x$. 

The homological index was  introduced by G\'omez-Mont in \cite{GM1} and it coincides  with the  GSV-index  introduced by  G\'omez-Mont,   Seade and  Verjovsky in  \cite{GSV}.    The concept of GSV-index has been extended to more general contexts, we refer to the works \cite{ a03, bss, Suw2, CoMa}.

Since $D$ is invariant by $\fol$, the interior multiplication $i_{v}$ also induces the complex of logarithmic differential forms
$$
0 \longrightarrow \Omega^n_{X,x}(\log\, \ddd) \stackrel{i_v}{\longrightarrow} \Omega^{n-1}_{X,x}(\log\, \ddd)\stackrel{i_v}{\longrightarrow} \cdots\stackrel{i_v}{\longrightarrow} \Omega^{1}_{X,x}(\log\, \ddd)\stackrel{i_v}{\longrightarrow} \OO_{n,x} \longrightarrow 0.
$$
If $x$ is an isolated singularity of $\fol$, the $i_v$-homology groups of the complex $\Omega^{\bullet}_{X,x}(\log\, \ddd)$ are finite-dimensional vector spaces (see \cite{Ale_p}). Thus, the Euler characteristic
$$
\chi (\Omega^{\bullet}_{X,x}(\log\, \ddd), i_v) = \displaystyle\sum_{i=0}^n(-1)^i \dim H_i(\Omega^\bullet_{X,x} (\log\, \ddd),i_v)
$$
\noindent of the complex of logarithmic differential forms is well defined.  Since this number does not depend on the local representative $v$ of the foliation $\fol$ at $p$, the {\it logarithmic index} of $\fol$  at the point $p$ is defined by
$$
\Log(\fol,\ddd,x):= \chi (\Omega^{\bullet}_{X,x}(\log\, \ddd), i_v).
$$
\noindent  Aleksandrov \cite{Ale_p} has proved the following formula that relates the logarithmic index and the residues  
\begin{eqnarray} \label{0003}
\Log(\fol, \ddd, x) = \mu(\fol,x) - \mbox{Ind}_{Hom} (v,  D, x),
\end{eqnarray}
where $\mbox{Ind}_{Hom}(v,  D, x)$ is the so-called  homological index of $v$ at $x$.  Recall that $\mbox{Ind}_{Hom}(v,  D, x)=\GSV(v,D,x)$ if $x$ is isolated, see \cite{GM1, GSV}.

\subsubsection{GSV-index and Camacho--Sad indices  on Surfaces}  \label{GSV-CS}

Let \( X \) be a complex compact surface and \( \fol \) a one-dimensional holomorphic foliation on \( X \). Let \( D \) be a reduced curve in \( X \). 
Consider \( \omega \in \mathrm{H}^0(X, \Omega_X^1 \otimes N_{\fol}) \), a rank-one Pfaffian system inducing \( \fol \). 
If \( D \) is invariant under \( \fol \), we say that \( \fol \) is {\it logarithmic along} \( D \).

Given a point \( x \in D \), let \( f = 0 \) be a local equation of \( D \) in a neighborhood \( U_{\alpha} \) of \( x \), and let \( \omega_{\alpha} \) be the holomorphic \(1\)-form inducing the foliation \( \fol \) on \( U_{\alpha} \). Since \( \fol \) is logarithmic along \( D \), it follows from \cite{Sai} that there are holomorphic functions \( g \) and \( \xi \) defined in a neighborhood of \( x \), which do not vanish identically simultaneously on \( D \), such that
\begin{eqnarray} \label{2p11}
g\,\frac{\omega_{\alpha}}{f} = \xi\, \frac{df}{f} + \eta,
\end{eqnarray}
where \( \eta \) is a suitable holomorphic \(1\)-form.
M. Brunella in  \cite{Bru97} showed that the GSV-index can be defined as follows.

\begin{defi}[Brunella \cite{Bru97}] Let \( \fol \) be a one-dimensional holomorphic foliation on a complex compact surface \( X \) and logarithmic along a reduced curve \( D \subset X \). Given \( x \in D \), we define
$$
\mathrm{GSV}(\fol, D, x) = \sum_i \mathrm{ord}_x \left(\frac{\xi}{g}\big|_{D_i}\right),
$$
where \( D_i \subset D \) are irreducible components of \( D \) and \( \mathrm{ord}_x \left(\frac{\xi}{g}\big|_{D_i}\right) \) denotes the order of vanishing of \( \frac{\xi}{g}\big|_{D_i} \) at \( x \).   
\end{defi}

 Given \( x \in D \), the Camacho--Sad index of $\fol$ at $x$ is defined as follows  
$$
\mathrm{CS}(\fol, D, x) = \sum_i \mathrm{ord}_x \left(-\frac{\eta}{\xi}\big|_{D_i}\right).
$$

\begin{teo}[Brunella \cite{Bru97}] Let \( \fol \) be a one-dimensional holomorphic foliation on a complex compact surface \( X \) and logarithmic along a reduced curve \( D \subset X \). Then
\begin{eqnarray}\nonumber
\sum_{x \in \Sing\left(\fol\right) \cap D} \mathrm{GSV}(\fol, D, x) = (N_{\fol}   - D) \cdot D\,.
\end{eqnarray}
\end{teo}

 \begin{teo}[Camacho--Sad \cite{Camacho-sad}] Let \( \fol \) be a one-dimensional holomorphic foliation on a complex compact surface \( X \) and logarithmic along a reduced curve \( D \subset X \). Then
\begin{eqnarray}\nonumber
\sum_{x \in \Sing\left(\fol\right) \cap D} \mathrm{CS}(\fol, D, x) = D^2\,.
\end{eqnarray}
\end{teo}

\subsection{The Relative $\check{\mbox{C}}$ech-de Rham cohomology} \label{s0005}
In this work, we will consider an appropriate modification in Chern-Weil theory so that it is adapted to the $\check{\mbox{C}}$ech-de Rham cohomology. For this purpose, next, we describe the relative $\check{\mbox{C}}$ech-de Rham cohomology. For more details, see, for example, \cite{bss} and \cite{Suw2}. 

Let $X$ be an $n$-dimensional complex manifold and $S$ a compact subset of $X$. We set the cover $\mathcal{U} = \{U_{0}, U_{1}\}$ of $X$, where $U_{0} = X - S$ and $U_{1}$ are  open neighborhood of $S$.  Considering $X$ as a real $2n$-dimensional oriented manifold, we will denote by $(A^{\bullet}(\mathcal{U}), d)$ the associated $\check{\mbox{C}}$ech-de Rham complex and by $A^{r}(\mathcal{U},U_{0})$ the kernel of the canonical projection $A^r(\mathcal{U}) \longrightarrow A^r(U_{0})$, $0 \leq r \leq 2n$. Since 
$$
A^{r}(\mathcal{U},U_{0}) = \{\xi = (\xi_0, \xi_1, \xi_{01}) \in A^r(\mathcal{U}) \,\,|\,\, \xi_0 = 0 \}, 
$$
\noindent we find that if $\xi \in A^{r}(\mathcal{U},U_{0})$, then $d\xi$ belongs to $A^{r+1}(\mathcal{U},U_{0})$. Thus, this gives rise to another complex, called the relative  $\check{\mbox{C}}$ech-de Rham complex.  The $r$-th relative $\check{\mbox{C}}$ech-de Rham cohomology of the pair $(\mathcal{U},U_{0})$ is defined as 
$$
H^r_d(\mathcal{U},U_{0}):= \mbox{Ker}\, d^r/\mbox{Im}\, d^{r-1}.
$$ 
\noindent By the five lema, there is an isomorphism $H^r_d(\mathcal{U},U_{0}) \cong H^r(X,X-S;\CC)$ (see \cite{Suw2}). Furthermore, if we assume that $U_1$ constitutes a regular neighborhood of $S$, then we have the Alexander duality (see \cite{brass01}).
$$
\mathcal{A}: H^{r}(X,X-S;\CC) \cong H^r_d(\mathcal{U},U_{0}) \longrightarrow H^{2n - r}(U_1,\CC)^{\ast} \cong H_{2n - r}(S,\CC)
$$
\noindent which defines the following commutative diagram

\begin{center}
\hspace{0.0 cm} \xymatrix{
H^{r}(X,X-S;\CC) \ar[r]^-{j^{\ast}} \ar[d]_{\mathcal{A}} &  H^{r}(X,\CC) \ar[d]^{\mathcal{P}}  \\
H_{2n - r}(S,\CC) \ar[r]^-{i_{\ast}} &  H_{2n-r}(X,\CC) \\
} 
\end{center}
\noindent where $i$ and $j$ denote, respectively, the inclusions $S \hookrightarrow X$ and $(X, \emptyset) \hookrightarrow (X, X - S)$ and 
$\mathcal{P}: H^{r}(X,\CC) \longrightarrow  H_{2n-r}(X,\CC) $ denote the Poincar\'e  duality.

Moreover, given $[\xi] = [(\xi_0,\xi_1,\xi_{01})] \in H^r_d(\mathcal{U},U_{0}) \cong H^r(X,X-S;\CC)$, we have that $\mathcal{A}([\xi]) \in H_{2n - r}(S,\CC)$ is represented by an $(2n - r)$ - cycle $C$ in $S$, such that for every closed $(2n-r)$-form $\eta$ on $U_1$  
$$
\displaystyle \int_C \eta =\int_{R_1}\xi_1\wedge \eta + \int_{R_{01}}\xi_{01}\wedge \eta,
$$
\noindent where $R_{1} \subset U_1$ be a real $2n$-dimensional manifold with $C^{\infty}$ boundary, compact, containing $S$ in its interior and $R_{01}=-\partial R_1$. In particular, if $r = 2n$, then $\mathcal{A}([\xi]) \in H_{0}(S,\CC) \cong \CC$ is a complex number given by
\begin{eqnarray}\label{eq03}
\mathcal{A}([\xi]) = \int_{R_1}\xi_1 + \int_{R_{01}}\xi_{01}.
\end{eqnarray}

\section{Proof of Theorem \ref{ttt1} and Corollary \ref{cor-princ-iso}}

Let $\fol$ be a one-dimensional holomorphic foliation on a complex manifold $X$, logarithmic along a free divisor $\ddd \subset X$ and $\varphi$ a homogeneous symmetric polynomial of degree $n$. Given a compact connected component $S_{\lambda}$ of $\Sing(\fol)$, let $U_{1\lambda}$ an open neighborhood of $S_{\lambda}$ in $X$ such that it does not intercept any other component of $\Sing(\fol)$. On $U_{0\lambda}:= U_{1\lambda} - S_{\lambda}$ we have an injective vector bundle homomorphism
$$
\sigma : T_{\fol}|_{U_{0\lambda}} \longrightarrow T_X(-\log \,\ddd)|_{U_{0\lambda}}
$$
\noindent  which associates each section $s$ of $T_{\fol}$, represented by a collection $(f_{\alpha})$ of holomorphic functions, to the vector field $v = f_{\alpha} v_{\alpha}$, which does not depend on the index $\alpha$ (by definition of $\fol$). Considering the quotient bundle $N_{\fol_0}:= T_X(- \log \,\ddd)|_{U_{0\lambda}}/ T_{\fol}|_{U_{0\lambda}}$, we have the following exact sequence of vector bundles on $U_{0\lambda}$ 
\begin{eqnarray}\label{seq0001}
0 \longrightarrow T_{\fol}|_{U_{0\lambda}} \stackrel{\sigma}{\longrightarrow} T_X(- \log \,\ddd)|_{U_{0\lambda}} \stackrel{\eta}{\longrightarrow} N_{\fol_0} \longrightarrow 0.
\end{eqnarray}
Furthermore, from the integrability of $\fol$ and since $T_X(-\log \,\ddd)$ is closed by the Lie bracket, we have an action of $T_{\fol}|_{U_{0\lambda}}$ on the vector bundle $N_{\fol_0}$ defined by 
\begin{center}
$\begin{array}{cccccc}
\theta : & \!  \Gamma (T_{\fol}|_{U_{0\lambda}} )\times \Gamma (N_{\fol_0})  & \! \longrightarrow & \!  \Gamma (N_{\fol_0}) \\
& \!  (s,\eta (\omega)) & \! \longmapsto & \!   \theta(s,\eta (\omega)) := \eta ([\sigma(s), \omega]).
\end{array}$
\end{center}
 \noindent Thus, considering a connection $\nabla_{\lambda}$ of type (1,0)  for $N_{\fol_0}$ satisfying 
$$
\nabla_{\lambda}(\eta(w))(s) = \eta([\sigma(s),w]),
$$
\noindent for all $w\in \Gamma (T_X(- \log \,\ddd)|_{U_{0\lambda}})$ and $s\in \Gamma (T_{\fol}|_{U_{0\lambda}})$, it follows from Bott vanishig theorem  \cite[Theorem 6.2.3]{bss}  that  $\varphi(\nabla_{\lambda}) \equiv 0$.

We set the covering $U_\lambda=\{U_{0\lambda},U_{1\lambda}\}$, where
$U_{0\lambda}:=U_\lambda\setminus S_\lambda$ and $U_{1\lambda}$ is an open neighborhood of $S_\lambda$ in $U_\lambda$.
 On the other hand, for each $k=0,1$, let $\nabla^{\bullet}_{k\lambda} = (\nabla'_{k\lambda},\nabla_{k\lambda})$ be a pair of connections for $T_{\fol}$ and $T_X(-\log \,\ddd)$, respectively, on $U_{k\lambda}$. Then the characteristic class $\varphi\Big(T_X(-\log \,\ddd) - T_{\fol}\Big)$ is represented by the cocycle 
$$\varphi(\nabla^{\bullet}_{\ast\lambda}) = \Big(\varphi(\nabla^{\bullet}_{0\lambda}), \varphi(\nabla^{\bullet}_{1\lambda}), \varphi(\nabla^{\bullet}_{0\lambda}, \nabla^{\bullet}_{1\lambda} )\Big)$$ 
\noindent in $\check{\mbox{C}}$ech-de Rham complex $A^{\bullet}(\mathcal{U_{\lambda}})$, for covering $\mathcal{U_{\lambda}} = \{U_{0\lambda}, U_{0\lambda}\}$. Now, choose $(\nabla'_{0\lambda},\nabla_{0\lambda})$ so that the triple $(\nabla'_{0\lambda},\nabla_{0\lambda}, \nabla_{\lambda})$ is compatible with the sequence (\ref{seq0001}) on $U_{0\lambda}$, i. e., such that the following relations occur 
$$(1\otimes\eta) \circ \nabla_{0\lambda}  = \nabla_{\lambda} \circ \eta \,\,\,\,\,\,\,\,\mbox{e} \,\,\,\,\,\,\,\,\,\,\,\,\,\,\,\,(1\otimes\sigma) \circ \nabla'_{0\lambda}  = \nabla_{0\lambda} \circ \sigma.
$$
Since the sequence (\ref{seq0001}) is exact, there is always triple of connections compatible with the sequence \cite[lema 4.17]{BB2}, . It follows from  \cite[Proposition 5.2.1]{bss}, that $\varphi(\nabla^\bullet_{0\lambda}) = \varphi(\nabla_{\lambda})$. Thus, we get that $\varphi(\nabla^\bullet_{0\lambda}) = 0$ and, consequently, the cocycle $\varphi(\nabla^{\bullet}_{\ast\lambda})$ is in relative  $\check{\mbox{C}}$ech-de Rham complex $A^{2n}(\mathcal{U_{\lambda}}, U_{0\lambda})$ and it defines a cohomology class 
$$
\varphi_{S_{\lambda}}(T_X(-\log \,\ddd) - T_{\fol}, \fol) \in  H^{2n}_d(\mathcal{U_{\lambda}}, U_{0\lambda}),
$$
\noindent that does not depend on the choice of the pair of connections $\nabla^{\bullet}_{1\lambda}$ or of triple of connections $(\nabla_{0\lambda}, \nabla'_{0\lambda}, \nabla_{\lambda})$.
Finally, the residue $\mbox{Res}^{\log}_{\varphi}\Big(\fol,D,  S_{\lambda}\Big)$ is defined as the image of class $\varphi_{S_{\lambda}}(T_X(-\log \,\ddd) - T_{\fol}, \fol)$ by the Alexander isomorphism. In other words,
$$
\mbox{Res}^{\log}_{\varphi}\Big(\fol,D,  S_{\lambda}\Big):= \mathcal{A}(\varphi_{S_{\lambda}}(T_X(-\log \,\ddd) - T_{\fol}, \fol)) \in H_0(S_{\lambda};\CC) \simeq \CC.
$$
and the item (i) of theorem is proved.

In order to prove item (ii), we will assume that $X$ is compact. Since $\Sing(\fol) = \displaystyle \bigcup_{\lambda} S_{\lambda}$, we have the following decomposition
\begin{eqnarray}\label{exxx}
H_0(\Sing(\fol);\CC) = \displaystyle\bigoplus_{\lambda}H_0(S_{\lambda};\CC),
\end{eqnarray}
\noindent such that for each $\lambda$ the residue $\mbox{Res}^{\log}_{\varphi}\Big(\fol,D,  S_{\lambda}\Big)$ belongs to the component $H_0(S_{\lambda};\CC)$. 
Thus, using the decomposition (\ref{exxx}) and the fact that $\Sing(\fol)$ is a compact subset of $X$ (since $X$ is compact), we obtain the following commutative diagram.  
\begin{center}
\hspace{0.0 cm} \xymatrix{
H^{2n}(X,X-\Sing(\fol);\CC) \ar[r]^-{j^{\ast}} \ar[d]_{\mathcal{A}} &  H^{2n}(X,\CC) \ar[d]^{\mathcal{P}}  \\
 \displaystyle H_0(\Sing(\fol);\CC) = \bigoplus_{\lambda}H_0(S_{\lambda};\CC)  \ar[r]^-{i_{\ast}}  \ar@{=}[dd]
&  H_0(X,\CC) \\
\\
 \displaystyle\bigoplus_{S_{\lambda} \not\subset D}H_0(S_{\lambda},\CC) \oplus  \displaystyle\bigoplus_{S_{\lambda} \subset D}H_0(S_{\lambda},\CC)   \ar@/_/[uur]_{i_*} & 
} 
\end{center} 
\noindent \noindent where $i$ and $j$ denote, respectively, the inclusions $\Sing(\fol) \hookrightarrow X$ and $(X, \emptyset) \hookrightarrow (X, X - \Sing(\fol))$. The map $\mathcal{A}: H^{2n}(X,X-\Sing(\fol);\CC) \longrightarrow H_0(\Sing(\fol),\CC)$ denote the Alexander isomorphism and $\mathcal{P}: H^{2n}(X,\CC) \longrightarrow  H_0(X,\CC)$ define the Poincar\'e duality. Hence by the commutativity of the diagram we get the desired formula in the homology $H_0(X,\CC)$:  
\begin{eqnarray} \nonumber
\int_{X} \varphi\Big(T_X(-\log\, \ddd) - T_{\fol}\Big) &= &  \sum_{\lambda} \mbox{Res}^{\log}_{\varphi}\Big(\fol,   \ddd  , S_{\lambda}\Big)\\ \nonumber
& = &\sum_{ S_{\lambda} \not\subset D } \mbox{BB}_{\varphi}\Big(\fol, S_{\lambda}\Big)+ \sum_{S_{\lambda}\subset  D   } \mbox{Res}^{\log}_{\varphi}\Big(\fol,   \ddd  , S_{\lambda}\Big),
\end{eqnarray}
\noindent where in the last equality, we have used the fact that the residue $\Res^{\log}_{\varphi}\Big(\fol,   \ddd  , S_{\lambda}\Big)$ coincides with the Baum-Bott residue $\mbox{BB}_{\varphi}\Big(\fol, S_{\lambda}\Big)$, when $S_{\lambda} \not\subset D$.

Finally, to prove item (iii), we consider $U \subset X$ as a neighborhood of the singular point $p \in \Sing(\fol) \cap \ddd$ and $\vartheta \in\Gamma ( T_X(-\log \,\ddd)|_{U})$ as a local representative of $\fol$.
The proof begins by demonstrating that the residue \( \mbox{Res}^{\log}_{\varphi}(\fol, D, p) \) coincides with the residue of the vector field \( \vartheta \) for the virtual vector bundle 
\[
T_X(-\log D) - T_{\fol}
\] 
at \( p \). This is followed by applying \cite[Theorem 6.2]{Suw2} to derive the desired formulas. In fact, let $U_0, U_1 \subset X$ be open sets defined by $U_0 := U - \{p\}$ and $U_1 := U$. Since $T_X(-\log \,\ddd)$ is closed under the Lie bracket, we have that $T_{X}(-\log\,\ddd)|_{U_0}$ becomes a holomorphic $\vartheta$-bundle by the action 
\begin{center}
$\begin{array}{cccc}
\alpha_\vartheta : & \! \Gamma (T_{X}(\log\,\ddd)|_{U_0} )  & \! \longrightarrow & \! \Gamma (T_{X}(\log\,\ddd)|_{U_0}) \\
& \! u & \! \longmapsto & \!  \alpha_{\vartheta}(u) = [\vartheta, u].
\end{array}$
\end{center}
In the open $U_0$ (restricting it if necessary), we can consider $T_{\fol}$ as the trivial line bundle $X \times \CC$. Thus, $T_{\fol}|_{U_0}$ becomes a holomorphic $\vartheta$-bundle by the action
\begin{center}
$\begin{array}{cccc}
\alpha'_{\vartheta} : & \! \Gamma (T_{\fol}|_{U_0} )  \cong \Gamma ({U_0})  & \! \longrightarrow & \! \Gamma (T_{\fol}|_{U_0}) \cong \Gamma ({U_0}) \\
& \! f & \! \longmapsto & \!  \alpha'_{\vartheta}(f) = \vartheta(f).
\end{array}$
\end{center}
\noindent Moreover, if we identify $T_{\fol}|_{U_0}$ with the subbundle $F_{\vartheta}$ of $T_{X}(\log\,\ddd)|_{U_0}$ spanned by $\vartheta$ on $U_0$ so that the trivialization $1$ given by $1(x) = 1$ corresponds to $\vartheta$, then $\alpha'_{\vartheta}$ corresponds to restriction of $\alpha_\vartheta$ to the sections of $\fol_{\vartheta}$.
Now, let $\nabla_0$ be a $\vartheta$ connection for $T_{X}(\log\,\ddd)$ on $U_0$ and $\nabla'_0$ and $\nabla$ connections for $T_{\fol}|_{U_0}$ and $N_{\fol_{0}}$, respectively, on $U_0$ obtained from $\nabla_0$ by restriction and going to the quotient. Then $\nabla$ is an $T_{\fol}|_{U_0}$ - connection and the triple $(\nabla'_0, \nabla_0, \nabla)$ is compatible with the sequence (\ref{seq0001}) on $U_0$. Hence, the residue of $\vartheta$ for $T_X(-\log\, \ddd) - T_{\fol}$ at the point $p$ coincides with $\mbox{Res}^{\log}_{\varphi}(\fol, D, p)$.  Now,  we can  use \cite[Theorem 6.2]{Suw2} to calculate the residue of $\vartheta$ for $T_X(-\log\, \ddd) - T_{\fol}$ at point $p$ : for the line bundle $T_{\fol}|_{U_0}$ we consider the frame defined by constant function $1(x) = 1$ and using the trivial relation
$$
\alpha'_{\vartheta}(1) =  0,
$$
we obtain the zero matrix $M^{1} = 0$ and hence, 
$$
\varphi(M^{1},M^{0}) = \varphi(M^{0}).
$$  
On the other hand, taking the frame $(\delta_1,\ldots, \delta_n)$ of $T_X(-\log\, \ddd)$ on $U_0$, obtained by restriction of the vector fields $\delta_1, \ldots, \delta_n$ that determine a system of $\mathcal{O}_{X,p}$-free basis for $T_{X,p}(-\log\, D)$ and considering  the expression of vector filed $\displaystyle \vartheta = \sum_{i=1}^{n} \vartheta_{i}\delta_{i}$, we have
$$
\alpha_{\vartheta}(\delta_j) =
[{\vartheta},\delta_{j}] = \sum_{i=1}^{n} [\vartheta_{i}\delta_{i}, \delta_{j}] = \sum_{i=1}^{n} m_{ij}\delta_{i}.$$
By applying the properties of the Lie bracket in
$$\displaystyle[\delta_{i}, \delta_{j}] = \sum_{k=1}^{n}\delta_{ij}^{k} \delta_{k}, \ \ \ \mbox{for} \ \ \ i < j,$$
\noindent  So, we set $M_{\log}(\vartheta):=(m_{ij})$ with entries given by \eqref{eq:Mlog-clean} and thus we obtain
$$\Res^{\log }_{\varphi}(\fol,D, p) = \Res_p\displaystyle\left[\begin{array}{cccc}\varphi\big(M(\vartheta)\big)  \\ v_1,\ldots,v_n \end{array}\right ].
$$
 Observe that if $i=j$ then $\delta_{ii}^{k} = 0$ for all $k$. Moreover, if $i > j$ then $\delta_{ij}^{k} = -\delta_{ji}^{k}$ by the anticommutativity property of the Lie bracket. So, we define $ M_{\log}(\vartheta):= M(\vartheta)$.

It follows from  \cite[Proposition B]{Faber} that  $T_{X,p}(-\log\, \ddd)$ is induced by an abelian Lie algebra of vector fields if and only if the germ $(D,p)$ is a normal crossing. 
Then $\delta_{ji}^{k}=0$, for all $i,j,k$. So $M(\vartheta) =  J_{\log} \vartheta$  and 
$$\Res^{\log }_{\varphi}(\fol,D, p) = \Res_p\displaystyle\left[\begin{array}{cccc}\varphi\big(J_{\log}(\vartheta)\big)  \\ v_1,\ldots,v_n \end{array}\right ].
$$
Now,  since  $\vartheta = \sum_{i=1}^nv_i\frac{\partial}{\partial z_i} = \sum_{i=1}^n\vartheta_i\delta_i$ we have $v_i = \sum_{j=1}^na_{ij}\vartheta_j,$ for all $i=1,\ldots, n$ and suitable holomorphic functions $a_{ij}$. In matrix form, we write $v = A \vartheta$
where $v=(v_i)_{n\times 1}$, $A=(a_{ij})_{n\times n}$ and $\vartheta=(\vartheta_i)_{n\times 1}$. Since
$$
J_{\log} \vartheta =[\delta_i(\vartheta_j)]= \left[\sum_{k=1}^{n} a_{ik} \dfrac{\partial}{\partial z_{k}}(\vartheta_j)\right]= A \cdot D\vartheta,
$$
we obtain the relation $$\det(J_{\log} \vartheta)=\displaystyle\det(A) \det\left(\frac{\partial\vartheta_i}{\partial z_j}\right).
$$
Therefore, we get
\begin{eqnarray}\nonumber
\Res^{\log }_{\det}(\fol,D, p) &=& \Res_p\displaystyle\left[\begin{array}{cccc}\det\big(J_{\log} \vartheta\big)  \\ v_1,\ldots,v_n \end{array}\right ]\\\nonumber &=& \Res_p\displaystyle\left[\begin{array}{cccc}\det(A) \det\left(\displaystyle\frac{\partial\vartheta_i}{\partial z_j}\right)  \\ v_1,\ldots,v_n \end{array}\right ]
\\\nonumber &=&\Res_p\displaystyle\left[\begin{array}{cccc}\det\left(\displaystyle\frac{\partial\vartheta_i}{\partial z_j}\right)  \\ \vartheta_1,\ldots,\vartheta_n \end{array}\right ],
\end{eqnarray} 
\noindent where, in the last equality, we use the transformation law of Grothendieck residues. Finally, we recall from \cite[Corollary 2]{Ale_p} that the logarithmic index is given by
\begin{eqnarray} \nonumber 
\Log(\F,D,p) =   \dim_{\CC} \frac{\OO_{X,p}}{\left\langle \vartheta_1, \ldots, \vartheta_n \right\rangle}=\Res_p\displaystyle\left[\begin{array}{cccc}\det\left(\displaystyle\frac{\partial\vartheta_i}{\partial z_j}\right)  \\ \vartheta_1,\ldots,\vartheta_n \end{array}\right ].   
\end{eqnarray}
This concludes the proof of Theorem \ref{ttt1} and Corollary \ref{cor-princ-iso} as desired.

\begin{remark} \label{rem:log-regular-vs-singular}
Let $(X,D)$ be a smooth surface with $D=\{xy=0\}$ in local coordinates, so that
$T_X(-\log D)$ is locally free with frame $e_1=x\partial_x$, $e_2=y\partial_y$.
A logarithmic vector field $v\in H^0(U,T_X(-\log D))$ can be written uniquely as
$v=f\,e_1+g\,e_2$ with $f,g\in\OO_{X,p}$.
Then $p$ is a  logarithmic singularity  of the induced logarithmic foliation iff
$f(p)=g(p)=0$.  In particular, if the vector field has a non-degenerate singularity, so    $J(v)=a\,x\partial_x+b\,y\partial_y $,  with $a,b\in \mathbb{C}^*$. Then 
$$\Res^{\log }_{\det}(v,D, p)= \frac{\det( J_{\log}(v))}{\det(v)}=\frac{0}{ab}=0,$$
since $\det( J_{\log}(v))=0.$
\end{remark}

 \begin{exe}\label{example-arrang}   We consider $X = \mathbb{P}^{3}$ with homogeneous coordinates $[z_{0}:z_{1}:z_{2}: z_{3}]$ and the divisor $\ddd = \{f = 0 \}$, where $f = z_{0} z_{1} z_{2} z_{3}.$ As usual, we denote $U_i=\{z_i\neq 0\}\subset \mathbb{P}^{3}$. 
On $U_{0}$  we have 
$\ddd|_{U_{0}} = \{x_{1}x_{2}x_{3}= 0 \}$  and   \ \ $T_{\mathbb{P}^{3}}(-\log \ddd)|_{U_{0}} = <\delta_{1},\delta_{2},\delta_{3}>$ with 

$$\delta_{1} = x_{1}\dfrac{\partial}{\partial x_{1}}, 
\delta_{2} = x_{2} \dfrac{\partial}{\partial x_{2}}, 
\delta_{3} = x_{3} \dfrac{\partial}{\partial x_{3}}.$$
Now, let $\fol$ be a one-dimension holomorphic foliation on $\mathbb{P}^{3}$ logarithmic along $D$ given in local chart by the following logarithmic vector field 
$$
\begin{array}{ll}
v|_{U_{0}}  & = (x_{1}-1)\delta_{1} + (x_{2}-1) \delta_{2}+(x_{3}-1) \delta_{3} \\ \\
 & =  x_{1}(x_{1}-1) \dfrac{\partial}{\partial x_{1}}  +x_{2}(x_{2}-1) \dfrac{\partial}{\partial x_{2}} +x_{3}(x_{3}-1)  \dfrac{\partial}{\partial x_{3}}  \\ \\
 & =  v_{1}\dfrac{\partial}{\partial x_{1}}+v_{2}\dfrac{\partial}{\partial x_{2}} +v_{3} \dfrac{\partial}{\partial x_{3}} ,
\end{array} $$
\noindent  whose  singular set is  
$\mbox{Sing}(v|_{U_{0}}) = \Big \{ p_{1} = (0,0,0); p_{2} = (0,0,1), p_{3} = (0,1,0), p_{4} = (1,0,0), p_{5} = (1,1,0), p_{6} = (1,0,1),p_{7} = (0,1,1), p_{8} = (1,1,1)  \Big\}.$ Note that only $P_{8} = (1,1,1)$ is not in $D$.  
Since $D$ is normal crossing, then  
$$M_{\log}(\vartheta)(x_{1},x_{2},x_{3})=J_{\log} \vartheta|_{U_{0}}(x_{1},x_{2},x_{3}) = \begin{pmatrix}     x_{1}  & 0 & 0   \\
                         0  & x_{2}  & 0  \\
                         0   &  0 & x_{3}
                                \end{pmatrix}.$$   
                                Thus $c_{1}(M|_{U_{0}}) =(x_{1}+x_{2}+x_{3}),$ and
$$\Res^{\log }_{c_{1}^3}(\fol,D, p_{i}) = \Res_{p_{i}} \displaystyle\left[\begin{array}{cccc}c_{1}^{3}\Big(M|_{U_{0}}\Big)dx_{1}\wedge dx_{2}\wedge dx_{3} \\ v_1v_{2}v_3 \end{array}\right ] = \dfrac{c_{1}^{3}(M|_{U_{0}})(p_{i})}{\det J v (p_{i})},$$ 
\noindent since  $\mbox{Sing}(v|_{U_{0}})$ is not degenerate. Then
$$\Res^{\log }_{c_{1}^3}(\fol,D, p_{i}) = \left\{\begin{array}{cl}
                            0,   & i=1  \\
                            1,  & i=2   \\
                            1,  & i=3   \\
                            1,  & i=4   \\
                           -8,  & i=5   \\
                           -8,  & i=6   \\
                           -8,  & i=7. 
\end{array}\right.$$
Since $p_{8} \notin D$ the residue associated at this point is the usual Baum-Bott Residue given by
$$\begin{array}{ll}
\Res_{c_{1}^3}(\fol, p_{8}) & =  \Res_{p_{8}} \displaystyle\left[\begin{array}{cccc}c_{1}^{3}\Big(Jv\Big)dx_{1}\wedge dx_{2}\wedge dx_{3} \\ v_1 v_{2} v_3 \end{array}\right ] \\ \\

&= \dfrac{c_{1}^{3}(Jv)(p_{8})}{\det J v (p_{8})} = 27.

\end{array}$$
Now, in   $U_{1}$, with  coordinates $(y_{1},y_{2},y_{3})$,  the logarithmic vector field is represented as  
$$ v|_{U_{1}} = y_{1}(y_{1}-1)\dfrac{\partial}{\partial y_{1}}+ y_{2}(y_{2}-1)\dfrac{\partial}{\partial y_{2}}+ y_{3}(y_{3}-1)\dfrac{\partial}{\partial y_{3}}
$$
\noindent  singular at $p_{9}= (0,0,0), p_{10} = (0,0,1), p_{11}= (0,1,0)$ and $p_{12} = (0,1,1)$ and similarly
$$\Res^{\log }_{c_{1}^3}(\fol,D, p_{i}) = \left\{\begin{array}{cl}
                            0,  & i=9   \\
                            1,  & i=10  \\
                            1,  & i=11  \\
                           -8,  & i= 12.
\end{array}\right.$$
In   $U_{2}$,  with coordinates  $(t_{1},t_{2},t_{3})$, we have the logarithmic vector field 
$$ v|_{U_{2}} = t_{1}(t_{1}-1)\dfrac{\partial}{\partial t_{1}}+ t_{2}(t_{2}-1)\dfrac{\partial}{\partial t_{2}}+ t_{3}(t_{3}-1)\dfrac{\partial}{\partial t_{3}}
$$
\noindent  singular at  $p_{13}= (0,0,0)$ and $p_{14}= (0,0,1)$ and
$$\Res^{\log }_{c_{1}^3}(\fol,D, p_{i}) = \left\{\begin{array}{cl}
                            0,  & i=13  \\
                            1,  & i=14.  
\end{array}\right.$$
Finally,  we represent the logarithmic vector field $v$ in  $U_{3}$, with coordinates $(w_{1},w_{2},w_{3})$, such as 
$$ v|_{U_{3}} = w_{1}(w_{1}-1)\dfrac{\partial}{\partial w_{1}}+ w_{2}(w_{2}-1)\dfrac{\partial}{\partial w_{2}}+ w_{3}(w_{3}-1)\dfrac{\partial}{\partial w_{3}}
$$

\noindent with an only  new   singularity $P_{15}= (0,0,0)$ and residue
$$\Res^{\log }_{c_{1}^3}(\fol,D, p_{15}) = 0.$$
Since $
\ddd = \cup_{i=0}^3\ddd_{i}
$, where $\ddd_{i} = \{ z_{i}=0 \} \subset \mathbb{P}^{3}$ and the short exact sequence
$$ 0 \longrightarrow \Omega^{1}_{\mathbb{P}^{3}} \longrightarrow \Omega^{1}_{\mathbb{P}^{3}}(\log \ddd) \longrightarrow \oplus_{i = 1}^{4}\mathcal{O}_{\ddd_{i}} \longrightarrow 0,$$
\noindent follows (see, \cite{DolKa}, Proposition 2.4) one has the Chern class of cotangent logarithmic bundle
$$c_{1}\Big(\Omega^{1}_{\mathbb{P}^{3}}(\log \ddd)\Big) = c_{1}(\Omega_{\mathbb{P}^{3}}^{1}) + c_{1}(\mathcal{O}_{\ddd_{0}})+ c_{1}(\mathcal{O}_{\ddd_{1}})+ c_{1}(\mathcal{O}_{\ddd_{2}})+ c_{1}(\mathcal{O}_{\ddd_{3}}). $$

Since   $c_{1}(\mathcal{O}_{\ddd_{i}})=h$ for $i=0,1,2,3$, and $c_{1}(\Omega_{\mathbb{P}^{3}}^{1}) = -c_{1}(T_{\mathbb{P}^{3}}) = -4h$, we obtain  
$$c_{1}(T_{\mathbb{P}^{3}}(-\log \ddd)) = -c_{1}\Big(\Omega^{1}_{\mathbb{P}^{3}}(\log \ddd)\Big) =0,$$
\noindent where $h$ is the hyperplane class.
Since the  foliation $\fol$ has degree 2,   its tangent bundle is  
$ T_{\fol} = \mathcal{O}(-1).$ In order to finish we check both sides of the Theorem \ref{ttt1}
\footnotesize
$$\begin{array}{ccl}

1=\displaystyle\int_{\mathbb{P}^{3}} c_{1}^{3}\Big(T_{\mathbb{P}^{3}}(-\log \ddd)- T_{\fol} \Big) & &   \\ \\
\displaystyle \sum_{p\in \Sing(\fol)   } \Res^{\log }_{c_{1}^3}(\fol,D, p) & = & 0+3\cdot 1-3\cdot  8+27+0+2\cdot 1-8+0+1+0=1.
\end{array}$$
Moreover, we observe that
$$
\sum_{p\in \Sing(\fol)   } \Res^{\log }_{c_{1}^3}(\fol,D, p)=1 >0
$$
and 
$$
4=\deg(D)<  \deg(\fol) +n=2+3.
$$
In accordance with Corollary \ref{Poincare-oddn}. 

\end{exe}

\subsection{Proof of  Corollary \ref{dim2-poinc}} 
We have the short exact sequence \eqref{eq:TX-logD-JD} and the flag of subsheaves
$T\fol\subset T_X(-\log D)\subset T_X$. Set
\[
N_{\fol}(\log D):=T_X(-\log D)/T\fol,\qquad
N_{\fol}:=T_X/T\fol.
\]
Then we get a commutative diagram with exact rows and columns:
\[
\begin{tikzcd}
&   & 0 \arrow[d] & 0 \arrow[d] \\
0 \arrow[r] & T\fol \arrow[r] \arrow[d,equal] & T_X(-\log D) \arrow[r] \arrow[d]
& N_{\fol}(\log D) \arrow[r] \arrow[d] & 0 \\
0 \arrow[r] & T\fol \arrow[r]   & T_X \arrow[r] \arrow[d]
& N_{\fol} \arrow[r] \arrow[d] & 0 \\
  &   & J_D(D) \arrow[r,equal] \arrow[d] & J_D(D) \arrow[r] \arrow[d] & 0 \\
&  & 0 & 0
\end{tikzcd}
\]
By the Snake lema, this yields a short exact sequence
\begin{equation}\label{eq:NFlogD-NF-JD}
0\longrightarrow N_{\fol}(\log D)\longrightarrow N_{\fol}
\longrightarrow J_D(D)\longrightarrow 0.
\end{equation}
From this we can extract an exact sequence 
 
\centerline{
\xymatrix{
0\ar[r]&  N_\mathscr{F}(\log\,  D)  \ar[r]&  \mathscr{I}_{\Sing(\fol)}\otimes N_{\fol} \ar[r] &  \ar[r]  \J_{D}(D)  &0  \\
}.
}
\noindent Then
$$
c_1(  N_\mathscr{F}(\log\,  D) )=N_{\fol}-D.
$$
Thus, by the classical Baum--Bott Theorem
$$
N_{\fol}^2=c_1^2( \mathscr{I}_{\Sing(\fol)}\otimes N_{\fol})=\sum_{x\in  \Sing(\fol)\cap(X\setminus D) } \BB_{c_1^2}(\fol, x)+\sum_{x\in  \Sing(\fol)\cap D }\BB_{c_1^2}(\fol, x).
$$
One the one hand, by Brunella  and  Camacho--Sad Theorems
\begin{equation*}\label{c_2-tg}
\begin{split}
c_1^2(  N_\mathscr{F}(\log\,  D)  &= (N_{\fol}-D)^2 \\
\\
&= N_{\fol}^2-2(N_{\fol}-D)\cdot D - D^2\\
\\
&=   \BB_{c_1^2}(\fol, X-D)+ \BB_{c_1^2} (\fol, D)-2\GSV(\fol,D)-\CS(\fol,D).
\end{split}
\end{equation*}
On the other hand, by Theorem \ref{ttt1}
\begin{equation*} 
\begin{split}
c_1^2(  N_\mathscr{F}(\log\,  D)  &= \BB_{c_1^2}(\fol, X-D)+ \BB_{c_1^2}^{\log}(\fol, D).
\end{split}
\end{equation*}
Therefore,
\begin{equation*} 
\begin{split}
\BB_{c_1^2}(\fol,D)- \BB_{c_1^2}^{\log}(\fol, D)&=2\GSV(\fol,D)+\CS(\fol,D). 
\end{split}
\end{equation*}
It follows form \cite[page 533]{Bru97} that 
\begin{equation*} 
\begin{split}
\BB_{c_1^2}(\fol,D)- \BB_{c_1^2}^{\log}(\fol, D)-\CS(\fol,D)&=2\GSV(\fol,D)\geq 0
\end{split}
\end{equation*}
is an obstruction for Poincar\'e problem.  Now, if $\fol$ is a generalize curve along $D$, then by  Brunella's result  \cite[Theorem]{Bru97}
$$
\GSV(\fol,D)=0
$$
and
$$
\BB_{c_1^2}(\fol,D)=\CS(\fol,D).
$$
Then
$$
\Res^{\log }_{c_1^2}(\fol, D)=0.$$

\subsection{Proof of Corollary \ref{Baum-Bott-sing}} 

Consider   a  functorial log  resolution  $\pi: (X,D) \to (Y, \emptyset)$ with excepcional divisor $D$, see  \cite[ Theorems 3.35, 3.34]{kollar07}.     
Since the singular locus of $X$ is invariant with respect to any automorphism, it follows from \cite[ Corollary 4.6]{GKK10}) that the twisted vector field $v\in H^0(Y,T_Y\otimes T_{\fol}^*)$ has a lift
$\tilde{v}\in H^0(X,T_X(-\log D)\otimes (T_{\pi^{-1}\fol})^*)$. Denoting by $\tilde{\fol}$ the foliation associated with the vector field $\tilde{v}$. 
 By Theorem \ref{ttt1} we have  
  \begin{multline*}
\int_{X} \varphi\big(T_{X}(-\log D) - (T_{\pi^{-1}\fol})^*) = \sum_{\pi^{-1}(x) \in \Sing(\pi^{-1} v) \cap (X \setminus D)} \BB_{\varphi}(\pi^{-1}v, \pi^{-1}(x))+ \\
\\
+ \sum_{S_{\lambda} \subset D \cap \Sing(\pi^{-1} v)} \Res^{\log}_{\varphi}\big(\pi^{-1} v, S_{\lambda}\big)= \\
\\
= \sum_{x \in \Sing(v) \cap Y_{\text{reg}}} \BB_{\varphi}(v, x) + \sum_{S_{\lambda} \subset D \cap \Sing(\pi^{-1} v)} \Res^{\log}_{\varphi}\big(\pi^{-1} v, S_{\lambda}\big).
\end{multline*}
\\
  \noindent  since  $\BB_{\varphi}(\pi^{-1}v, \pi^{-1}(x))= \BB_{\varphi}( v, x)$, for all $x\in Y_{\reg} $.  If,   \( \pi^*\fol \) is induced by a global vector field and \( \Sing(\fol) \subset \Sing(Y) \), then 
\footnotesize
 $$\chi (Y)-\chi (\Sing(Y))=\chi (X-D)= \int_{X} c_n\big(T_{X}(-\log D) = \sum_{ x \in D\cap \Sing(\pi^{-1}v)} \Res_x\displaystyle\left[\begin{array}{cccc}\varphi(J_{\log}(\pi^{-1}v))  \\ v_1,\ldots,v_n \end{array}\right ].
 $$

\begin{exe} 
Let $\wP^2_k:=\wP(1,1,k)$, $k\ge1$, and let $\pi:(\Sigma_k,D)\to\wP^2_k$ be the minimal good resolution,
where $\Sigma_k=\wP(\OO_{\wP^1}\oplus \OO_{\wP^1}(k))$ and $D\simeq\wP^1$ is the exceptional section with $D^2=-k$.
Pick $a_0,a_1,a_2\in\C^*$ with $a_0\neq a_1$, $ka_0\neq a_2$, $ka_1\neq a_2$, and let $\fol$ be the
one--dimensional foliation on $\wP^2_k$ induced by
$$\xi_a=a_0z_0\partial_{z_0}+a_1z_1\partial_{z_1}+a_2z_2\partial_{z_2}.$$
Then $\fol$ has isolated singularities and its pullback foliation $\widetilde\fol$ on $\Sigma_k$
has two non-degenerate singularities in $\Sigma_k\setminus D$ and (possibly) finitely many singularities on $D$,
all of them non-degenerate.
Since $D$ is smooth, then   the logarithmic indices along $D$ vanish at such points,
hence $\Ind_{\log}(\widetilde\fol,D,p)=0$ for every $p\in\Sing(\widetilde\fol)\cap D$.
Therefore,  
\[
\int_{\Sigma_k} c_2\!\bigl(T_{\Sigma_k}(-\log D)\bigr)
=\sum_{x\in \Sing(\widetilde\fol)\cap(\Sigma_k\setminus D)}\mu_x(\widetilde\fol)=2.
\]
On the other hand, using $c_2(T_{\Sigma_k})=4$, $K_{\Sigma_k}=-2D-(k+2)L$, $D^2=-k$ and $D\cdot L=1$, one gets
\[
\int_{\Sigma_k} c_2\!\bigl(T_{\Sigma_k}(-\log D)\bigr)
=c_2(T_{\Sigma_k})+K_{\Sigma_k}\cdot D + D^2=4+(k-2)-k=2,
\]
confirming the formula in this explicit resolution.
\end{exe}

 \section{Proof of Corollary \ref{mu-PHGB}}
\begin{proof}
By Theorem \ref{teo_1.1} applied to $(\widehat X,\widehat D)$ and the logarithmic foliation $\widehat\fol$,
\begin{equation}\label{eq:BBlog-short}
\footnotesize
\int_{\widehat X} c_n\!\bigl(T_{\widehat X}(-\log \widehat D)-T_{\widehat\fol}\bigr)
=
\sum_{x\in \Sing(\widehat\fol)\cap(\widehat X\setminus \widehat D)}\mu_x(\widehat\fol)
+
\sum_{p\in \Sing(\widehat\fol)\cap \widehat D}\Log(\widehat\fol,\widehat D,p).
\end{equation}
Since $\widehat{v}$ is a global section, $c(T_{\widehat\fol})=1$, hence
$c_n\!\bigl(T_{\widehat X}(-\log \widehat D)-T_{\widehat\fol}\bigr)=c_n\!\bigl(T_{\widehat X}(-\log \widehat D)\bigr)$.
By the logarithmic Poincar\'e--Hopf theorem on the complement  and
$\widehat X\setminus \widehat D \simeq X\setminus D$, we have
\begin{equation}\label{eq:PH-short}
\sum_{x\in \Sing(\widehat\fol)\cap(\widehat X\setminus \widehat D)}\mu_x(\widehat \fol )=\chi(X\setminus D).
\end{equation}
Substituting into \eqref{eq:BBlog-short} yields
\begin{equation}\label{eq:indices-int-minus-chi-short}
\sum_{p\in \Sing(\widehat\fol)\cap \widehat D}\Log(\widehat\fol,\widehat D,p)
=
\int_{\widehat X} c_n\!\bigl(T_{\widehat X}(-\log \widehat D)\bigr)-\chi(X\setminus D).
\end{equation}

Let $\widehat Y:=\widehat D_{\mathrm{sing,sch}}$ and $\mu_L(\widehat Y)$ be the $\mu$--class of
\cite{Aluffi1}. From \cite[Corollary~1.5]{FerreiraLourenco2024},
\begin{equation}\label{eq:FL-short}
\chi(X\setminus D)
=
\int_{\widehat X} c_n\!\bigl(T_{\widehat X}(-\log \widehat D)\bigr)
+
(-1)^{n+1}\int_{\widehat X}\mu_LC(\widehat Y).
\end{equation}
Inserting \eqref{eq:FL-short} into \eqref{eq:indices-int-minus-chi-short}; the logarithmic Chern integrals cancel and we obtain
\eqref{eq:log-index-mu-class-short} as desired.
\end{proof}

\section{Proof of Corollary \ref{GZL}}

\begin{proof}
Let $v\in H^0(Y,T_Y\otimes\mathscr L)$ be a twisted vector field inducing the foliation $\fol$ on $Y$, and assume
$\Sing(\fol)\subset \Sing(Y)$. Fix a birational morphism $\pi:X\to Y$ obtained as a composition of a log resolution of $Y$
and a resolution of singularities of the induced foliation, so that $X$ is a smooth compact complex surface, the reduced
$\pi$-exceptional divisor $D\subset X$ is a simple normal crossing curve and coincides with the exceptional locus of $\pi$,
and the induced foliation $\G$ on $X$ is logarithmic along $D$, with $\Sing(\G)\subset D$. By hypothesis $\fol$ is a
generalized curve, hence so is $\G$ along $D$.

Pulling back $v$ to $X$ we obtain  a logarithmic twisted vector
field $$v_X\in H^0\bigl(X,\,T_X(-\log D)\otimes\pi^*\mathscr L\bigr).$$ The section $v_X$ may vanish divisorially along $D$,
but the saturated rank-one subsheaf of $T_X(-\log D)$ generated by $v_X$ defines $\G$; equivalently there exists an effective
divisor $Z$ supported on $D$ such that the cotangent line bundle of $\G$ is
$ \pi^*\mathscr L\otimes\mathcal O_X(Z)$.

Consider the logarithmic normal bundle $N_\G(\log D):=T_X(-\log D)/T_\G$. Since $X$ is smooth and $\G$ is logarithmic along
$D$, Corollary~\ref{dim2-poinc} applies to $(X,\G,D)$ and, because $\G$ is a generalized curve along $D$, yields
$$c_1\bigl(N_\G(\log D)\bigr)=\Res^{\log}_{c_1^2}(\G,D)=0.$$ On the other hand, using
$c_1\bigl(T_X(-\log D)\bigr)=-(K_X+D)$ and $c_1(T_\G^*)=c_1(\pi^*\mathscr L)+Z$, we compute
$c_1\bigl(N_\G(\log D)\bigr)=-(K_X+D)+\pi^*c_1(\mathscr L)+Z$. Since $T_Y$ is locally free, $\omega_Y$ is invertible and
$K_Y$ is Cartier; writing $K_X=\pi^*K_Y+\sum_i a_iD_i$ with $a_i\in\mathbb Z$ and $Z=\sum_i z_iD_i$ with $z_i\ge 0$, we obtain
$$c_1\bigl(N_\G(\log D)\bigr)=\pi^*(c_1(\mathscr L)-K_Y)+\sum_i(z_i-1-a_i)D_i.$$

Set $H:=\pi^*(c_1(\mathscr L)-K_Y)$ and $E:=\sum_i(z_i-1-a_i)D_i$, so that $$c_1\bigl(N_\G(\log D)\bigr)=H+E.$$ Each $D_i$ is
$\pi$-exceptional, hence $H\cdot D_i=0$, and the equality $c_1\bigl(N_\G(\log D)\bigr)=0$ implies
$0=(H+E)\cdot D_i=E\cdot D_i$ for all $i$. Since $D$ is the reduced exceptional divisor of a resolution of a normal surface,
the intersection matrix $(D_i\cdot D_j)$ is negative definite, hence invertible; the relations $E\cdot D_i=0$ for all
$i$ force $E=0$. Consequently $z_i-1-a_i=0$ for every $i$, so $a_i=z_i-1\ge -1$. Thus $Y$ has log canonical singularities.
Finally, the Zariski--Lipman conjecture is known for log canonical surfaces \cite{Dru14}. 
\end{proof}

\begin{exe}\label{ex:boundary-3-points}
Let $U\subset \mathbb{C}^3_{x,y,z}$ be a small polydisc around $0$, and set
$$
D:=\{xyz=0\}\cap U,
$$
so $D$ is SNC and $T_U(-\log D)$ is locally free with frame
$(x\partial_x,y\partial_y,z\partial_z)$.
Let $\pi:\widehat U\to U$ be the blow-up of the origin and denote by
$E\simeq \mathbb{P}^2$ the exceptional divisor. Let
$\widehat D:=\pi^{-1}(D)_{\mathrm{red}}=E\cup \widetilde D$, where
$\widetilde D$ is the reduced strict transform of $D$. Then $\widehat D$ is SNC.
In the standard affine chart $U_x=\{x\neq 0\}$ on $\widehat U$ with coordinates
$$
x,\qquad u:=y/x,\qquad v:=z/x,
$$
one has $\widehat D\cap U_x=\{xuv=0\}$ and $T_{\widehat U}(-\log\widehat D)$ has local frame
$$
x\partial_x,\qquad u\partial_u,\qquad v\partial_v.
$$
Similarly, in the charts $U_y$ and $U_z$ one has analogous coordinates
$(y,\,u':=x/y,\,v':=z/y)$ and $(z,\,u'':=x/z,\,v'':=y/z)$, with $\widehat D$ given by
$\{y u' v'=0\}$ and $\{z u'' v''=0\}$, respectively.

Define a logarithmic vector field $\widehat v$ by the following local expressions:
\begin{equation}\label{eq:vhat-Ux}
\widehat v|_{U_x}:=x\cdot(x\partial_x)+u\cdot(u\partial_u)+v\cdot(v\partial_v),
\end{equation}
\begin{equation}\label{eq:vhat-Uy}
\widehat v|_{U_y}:=y\cdot(y\partial_y)+u'\cdot(u'\partial_{u'})+v'\cdot(v'\partial_{v'}),
\end{equation}
\begin{equation}\label{eq:vhat-Uz}
\widehat v|_{U_z}:=z\cdot(z\partial_z)+u''\cdot(u''\partial_{u''})+v''\cdot(v''\partial_{v''}).
\end{equation}
These local formulas agree on overlaps (they come from the infinitesimal generator of the diagonal
$\mathbb{C}^*$-action), hence define a global section $\widehat v$ of $T_{\widehat U}(-\log\widehat D)$.
In each chart, $\widehat v$ vanishes exactly at the origin of that chart:
$$
(x,u,v)=(0,0,0)\in U_x,\qquad (y,u',v')=(0,0,0)\in U_y,\qquad (z,u'',v'')=(0,0,0)\in U_z.
$$
These three points lie on $E$ and correspond to the three coordinate points
$$
[1:0:0],\quad [0:1:0],\quad [0:0:1]\in E\simeq \mathbb{P}^2.
$$
There are no other zeros of $\widehat v$ on $\widehat U$.

In SNC coordinates, the logarithmic index at an isolated zero is the length of the local algebra cut out by the
coefficients in the logarithmic frame. For instance, in $U_x$ we have coefficients $(x,u,v)$, hence
$$
\Log(\widehat v,\widehat D,p_x)
=\dim_{\mathbb{C}}\frac{\mathbb{C}\{x,u,v\}}{(x,u,v)}=1,
$$
and similarly $\Log(\widehat v,\widehat D,p_y)=1$ and
$\mathrm{Ind}_{\log}(\widehat v,\widehat D,p_z)=1$. Therefore
\begin{equation}\label{eq:sum-indices-3}
\sum_{p\in \mathrm{Sing}(\widehat v)\cap \widehat D}
\Log(\widehat v,\widehat D,p)=3.
\end{equation}

Let $\widehat Y:=\widehat D_{\mathrm{sing,sch}}$ and let $L:=\mathcal{O}_{\widehat U}(\widehat D)$.
Formally applying Corollary~\ref{mu-PHGB} with $n=3$ gives
$$
\sum_{p\in \mathrm{Sing}(\widehat v)\cap \widehat D}\Log(\widehat v,\widehat D,p)
= - \int_{\widehat U}\mu_L(\widehat Y),
$$
hence, by \eqref{eq:sum-indices-3},
$$
\int_{\widehat U}\mu_L(\widehat Y)=-3.
$$
\end{exe}

\noindent
If one considers a local DT-type situation in which the relevant moduli space is locally modeled as a critical locus
$M=\mathrm{Crit}(W)\subset \mathbb{C}^3$, then the local contribution to the Behrend weighted Euler characteristic is the
Behrend value $\nu_M(0)$. \cite{Behrend2009}
For an isolated hypersurface critical point one has
$$\nu_M(0)=\mu(W,0)=\dim_{\mathbb{C}}\mathbb{C}\{x,y,z\}/(W_x,W_y,W_z),$$ 

\noindent see \cite{Behrend2009,Milnor}. 
Thus, in concrete problems where the local Behrend contribution equals $3$, the computation
\eqref{eq:sum-indices-3} provides a boundary-local expression producing the same integer.
\section{Residue-valued constructible functions}\label{sec:residue-constructible}

Let $X$ be a complex manifold. Write $CF(X)$ for the abelian group of $\Z$-valued
constructible functions on $X$, i.e.\ finite $\Z$-linear combinations of indicator
functions of constructible subsets.
For any commutative ring $R$ set
\begin{equation}\label{eq:CF-coeff}
CF(X;R)\;:=\;CF(X)\otimes_{\Z}R.
\end{equation}
Elements of $CF(X;R)$ may be viewed as $R$-valued constructible functions on $X$.
In particular, $CF(X;\C)$ is the natural receptacle for constructible functions with
complex values. If $c_*:CF(X)\to H_*^{BM}(X;\Z)$ denotes MacPherson's transformation,
then it extends by $R$-linearity to $CF(X;R)$; see \cite{MacPherson,SchuermannBook}.

Fix a simple normal crossings divisor $D\subset X$ and a one-dimensional logarithmic
foliation $\F$ along $D$. Assume that $\Sing(\F)$ is finite.
Let $\varphi$ be a homogeneous symmetric polynomial of the degree required in our
logarithmic Baum--Bott residue formalism.

For every $x\in\Sing(\F)$, our definition of the logarithmic Grothendieck--Baum--Bott
residue produces a complex number
\begin{equation}\label{eq:local-residue}
\Res^{\log}_{\varphi}(\F,D;x)\in\C.
\end{equation}
Define the \emph{residue-valued constructible function} $\varphi_{\F,D}$ on $X$ by
\begin{equation}\label{eq:varphiFD-def}
\varphi_{\F,D}(y)\;:=\;
\begin{cases}
\Res^{\log}_{\varphi}(\F,D;y), & y\in\Sing(\F),\\
0, & y\notin\Sing(\F).
\end{cases}
\end{equation}
Equivalently,
\begin{equation}\label{eq:varphiFD-sum}
\varphi_{\F,D}
\;=\;
\sum_{x\in\Sing(\F)} \Res^{\log}_{\varphi}(\F,D;x)\cdot \mathbf 1_{\{x\}}
\qquad\in\qquad CF(X;\C).
\end{equation}

Moreover, the assignment $\varphi\mapsto \varphi_{\F,D}$ is $\C$-linear in the polynomial:
if $\varphi=\sum_{\nu} c_{\nu}\,\varphi_{\nu}$ is an expansion in any fixed basis of
homogeneous symmetric polynomials, then
\begin{equation}\label{eq:linearity}
\varphi_{\F,D}\;=\;\sum_{\nu} c_{\nu}\,(\varphi_{\nu})_{\F,D}
\qquad\text{in }CF(X;\C).
\end{equation}
Let $i_x:\{x\}\hookrightarrow X$ be the inclusion. Set
\begin{equation}\label{eq:G}
\mathcal G_{\varphi}
\;:=\;
\bigoplus_{x\in\Sing(\F)} (i_x)_*(\C),
\end{equation}
viewed as a constructible sheaf supported on $\Sing(\F)$ and placed in degree $0$.
Define an endomorphism $\Phi_{\varphi}:\mathcal G_{\varphi}\to\mathcal G_{\varphi}$
by prescribing its stalk at $x$ to be multiplication by the scalar
$\Res^{\log}_{\varphi}(\F,D;x)$:
\begin{equation}\label{eq:Phi}
(\Phi_{\varphi})_x:\C\to\C,
\qquad
(\Phi_{\varphi})_x(\lambda)=\Res^{\log}_{\varphi}(\F,D;x)\cdot\lambda.
\end{equation}
 
The trace function associated with $(\mathcal G_{\varphi},\Phi_{\varphi})$ coincides
with $\varphi_{\F,D}$, i.e.
\begin{equation}\label{eq:trace-identification}
\varphi_{\mathcal G_{\varphi},\Phi_{\varphi}}=\varphi_{\F,D}
\qquad\text{in }CF(X;\C).
\end{equation}
Since $\mathcal G_{\varphi}$ is concentrated in degree $0$ and supported on $\Sing(\F)$,
for each $x\in\Sing(\F)$ we have $H^0((\mathcal G_{\varphi})_x)=\C$ and
$H^i((\mathcal G_{\varphi})_x)=0$ for $i\neq 0$. Hence the definition of the trace
function yields
\begin{equation*}
\varphi_{\mathcal G_{\varphi},\Phi_{\varphi}}(x)
=
\mathrm{Tr}\bigl((\Phi_{\varphi})_x\mid \C\bigr)
=
\Res^{\log}_{\varphi}(\F,D;x),
\end{equation*}
and $\varphi_{\mathcal G_{\varphi},\Phi_{\varphi}}(y)=0$ for $y\notin\Sing(\F)$.
This is exactly \eqref{eq:varphiFD-def}.

Characteristic cycles $CC(\mathcal H)$ carry \emph{integral} multiplicities; consequently,
invariants extracted solely from $CC(\mathcal H)$ are typically integer-valued.
Complex coefficients enter naturally once one enriches the coefficient system by an
endomorphism and passes to trace functions, as in \eqref{eq:trace-identification};
compare the microlocal Lefschetz philosophy for pairs $(\mathcal H,\Psi)$ in
\cite[Ch.\ 9]{KashiwaraSchapira}.

\end{document}